\documentclass{report}

\usepackage{amsfonts,amsmath,eufrak}
\usepackage{epsfig}
\usepackage{hyperref}
\usepackage{color}

\usepackage{graphicx}
 
\newtheorem{lemma}{Lemma}[chapter]
\newtheorem{theorem}[lemma]{Theorem}
\newtheorem{definition}[lemma]{Definition}

\newcommand{\D}[1]{\ensuremath{\mathbf{D}_{#1}}}
\newcommand{\K}[1]{\ensuremath{\mathbf{k}_{#1}}}
\newcommand{\LL}{\ensuremath{\widetilde{\K{p/q}}}}

\newcommand{\C}{\ensuremath{\mathbb{C }}}
\newcommand{\Z}{\mathbb{Z }}
\newcommand{\R}{\mathbb{R}}

\newcommand{\puq}{\ensuremath{\mathcal{PUQ}}}
\newcommand{\Q}{\ensuremath{\mathbb{H}}}
\newcommand{\quat}[4]{#1 + #2 {\bf i} + #3 {\bf j} + #4 {\bf k}}
\newcommand{\qvec}[1]{{\bf #1}}

\begin{document}
\pagenumbering{roman}
\begin{center}
{\Large \bf Great Circle Links in the Three-Sphere}
\vskip 0.2in
\centerline{By}
\vskip 0.1in
\centerline{\large  Genevieve Walsh}
%\centerline{\large B.A. (Oberlin College) 1994}
%\centerline{\large M.S. (Auburn University) 1997}
%\vskip 0.3in
%\centerline{\large  DISSERTATION}
\vskip 0.1in
\centerline{\large Submitted in partial satisfaction of the requirements for the
degree of}
\vskip 0.1in
\centerline{\large  DOCTOR OF PHILOSOPHY}
\vskip 0.1in
\centerline{\large in}
\vskip 0.1in
\centerline{\large MATHEMATICS}
\vskip 0.2in
\centerline{\large in the}
\vskip 0.1in
\centerline{\large  OFFICE OF GRADUATE STUDIES}
\vskip 0.1in
\centerline{\large of the}
\vskip 0.1in
\centerline{\large  UNIVERSITY OF CALIFORNIA}
\vskip 0.1in
\centerline{\large  DAVIS}
\end{center}
\vskip 0.1in
%{\large Approved:}
\begin{center}
%\centerline{\underbar{\hskip 2.5in}}
%\vskip 0.15in
%\centerline{\underbar{\hskip 2.5in}}
%\vskip 0.15in
%\centerline{\underbar{\hskip 2.5in}}
%\vskip 0.2in
%\centerline{\large Committee in Charge}
%\vskip 0.2in
\centerline{\large June 2003}
\end{center}
\newpage
\large
%\listoffigures
\tableofcontents

%%%%%%%%%%%%%%%%%%%%%%%%%%%%%%%%%%%%%%%%%%%%%%%%%%%%%%%%%%%%%%%%%%
%%%%%%%%%%%%%%%%%%%%%%%%%%%%%%%%%%%%%%%%%%%%%%%%%%%%%%%%%%%%%%%%%%

\newpage
\large
{\Large \bf ACKNOWLEDGMENTS} \\

I have been very fortunate to be in the creative and supportive environment of the UC Davis math department. I would particularly like to thank my advisor, Bill Thurston, for suggesting the problems in this dissertation and for many very interesting and helpful mathematical conversations. Abby Thompson and Joel Hass are also due special thanks for their advice, mathematical and otherwise.  Finally, I would like to thank my husband, Glen Nuckolls,  for his encouragement and understanding.  

This research was supported in part by  NSF VIGRE grant DMS-0135345.

%%%%%%%%%%%%%%%%%%%%%%%%%%%%%%%%%%%%%%%%%%%%%%%%%%%%%%%%%%%%%%%%%%%%%%%%%%%%%%%%%%%%%%%%%%%%%%%%%%%%%%%%%%%%%%%%%%%%%%%%%%%%%%%%%%%%%%%%%%%%%%%%%%%%%%%%%%

\newpage 
\begin{center}
\underline{\bf \Large Abstract} 
\end{center}
We investigate great circle links in the three-sphere, the class of links where each component is a great circle. 
Using the geometry of their complements,  we classify such  links up to five components.  Every two-bridge knot complement is covered by a great circle link complement.  We show that for roughly half of all two-bridge knots,  this cover contains a closed incompressible surface. Infinitely many fillings of the two-bridge knot lift to fillings of the great circle link  where the incompressibility of this surface is preserved.   It has been shown in \cite{hatcherthurston} that only finitely many fillings of a two-bridge knot produce Haken manifolds.  A closed manifold is said to be virtually Haken if has a finite cover that is Haken, but is not itself Haken. 
Thus, infinitely many fillings of an infinite class of  two-bridge knot complements are explicitly shown to be virtually Haken. 
\medskip

\large

%%%%%%%%%%%%%%%%%%%%%%%%%%%%%%%%%%%%%%%%%%%%%%%%%%%%%%%%%%%%%%%%

\newpage
\pagestyle{myheadings} 
\chapter{Introduction}
\pagenumbering{arabic}
\markright{  \rm \normalsize CHAPTER 1. \hspace{0.5cm}
 Introduction}
\thispagestyle{myheadings}
Knots and links are integral in the study of three-manifolds. For example, we know that every closed  three-manifold can be written as Dehn surgery on a link in $S^3$ \cite{Lickorish}.  Here we study a subclass of links, great circle links, where every component of the link is a geodesic in the standard metric on $S^3$.   This class of links has very natural connections to other areas of mathematics.  We hope to show that it is also very beautiful in its own right. 

\section{Geodesics in the three-sphere}
First we claim that every geodesic $g$ in $S^3$  is the intersection of a plane in $\R^4$ with the unit three-sphere.   
 Let $\qvec{x},\qvec{y}$ be two points on $g$, close enough together so that a geodesic between them is unique. Then the reflection of $\R^4$ though any 3-plane that contains $\qvec{x}, \qvec{y}$ and the origin in $\R^4$ induces an isometry on $S^3$ that fixes \qvec{x} and \qvec{y}.  Therefore this isometry fixes $g$.  Hence $g$ is in every such 3-plane, and in the unique two-plane through \qvec{x},\qvec{y} and the origin in $\R^4$.  The intersection of this 2-plane and $S^3$ is a curve on $S^3$ which we call a great circle.  This and the general theory of the geometry of the three-sphere can be found in the excellent survey,  \cite{Scottsurvey}.    

The author knows of no previous efforts to classify great circle links using their complements in $S^3$.  However,  great circle links have been studied fairly extensively in other contexts.  For example, a great circle in $S^3$ is invariant under the action of $\phi: S^3 \rightarrow S^3$ by $\phi(\qvec{x}) = -\qvec{x}$.  Therefore, a great circle link of $n$ components corresponds to a configuration of $n$ lines in $\mathbb{RP}^3$. Configurations of lines in $\mathbb{RP}^3$ have been classified, up to 6 components, by Viro and Drobotukhina in \cite{RP3}, using totally different, combinatorial  methods.  

 In addition, since every geodesic is naturally identified with a plane through the origin, a link of $p$ great circles corresponds to $p$ points in the Grassmanian of two-planes in $\R^4$, $G(4,2)$.    In \cite{conwaypackings},  Conway, Hardin, and Sloan investigate  packings of $G(m,n)$.  For $G(4,2)$ they find a configuration of six points that corresponds to a great circle link with remarkable symmetry properties.  Namely,  every two circles can be taken to any other two through an orientation-preserving symmetry.  
 
In chapter 2, we classify great circle links up to five components. In general we take the geometric point of view given in \cite{Thurstonbook} as the set-up for the classification. 

\section{Connections with rational knots}
We also relate great circle links to rational knots. A two-bridge knot  is a knot that can be isotoped to have two minima and two maxima, and no less, with respect to the standard height function on $R^3$.  By Schubert \cite{schubertknot} every two-bridge knot can be put in  a standard position,  consisting of two arcs of slope $p/q$ connecting punctures on a  four-punctured sphere.   In this way a two-bridge knot is associated with a (non-unique) rational number.   Thus we will use the term ``rational knot'' interchangeably with ``two-bridge knot''.   A sketch of Schubert's proof, from a slightly different perspective, is given in section 3.4.

An {\it incompressible surface} $S \subset M$ is an embedded surface $S \neq S^2$ in a three-manifold $M$ where  no non-trivial  embedded curve on $S$ bounds an embedded disk $D$ in $M$, with $D \cap S = \partial D$.   An orientable  irreducible three-manifold  that contains an  incompressible surface is called {\it Haken}.  Haken manifolds are a very well-understood class of three-manifolds.  Using Haken's notion of a hierarchy, Waldhausen \cite{waldhausenhaken} showed that Haken manifolds are characterized by their fundamental groups. Specifically, let $M$ and $N$ be irreducible 3-manifolds where any boundary components are incompressible and $M$ is Haken.  If  $\phi: \pi_1(N) \rightarrow \pi_1(M)$ is an isomorphism that respects the peripheral structure, then $\phi$ is induced by a homeomorphism $f: N \rightarrow M$.  ($\phi$ respects the peripheral structure if for every boundary component $F \subset N$, the subgroup of $\pi_1(N)$  corresponding to the inclusion of $\pi_1(F)$ is mapped to a subgroup that is conjugate to the inclusion of the fundamental group of a boundary component of $M$. )   W. Thurston also showed that Haken manifolds satisfy his geometrization conjecture.  Namely, Haken manifolds can be cut along incompressible tori so that each piece is either a Seifert fibered space or a hyperbolic 3-manifold.  A proof of this can be found in the book detailing the proof on the Smith Conjecture, \cite{Smithconj.}.   

However, not all manifolds are Haken. In \cite{hatcherthurston}, Hatcher and Thurston showed that all but finitely many fillings of rational knot complements are non-Haken.  We say that a manifold is {\it virtually Haken} if it is not itself Haken, but has a finite cover that is Haken. 
If  manifold  is virtually Haken then it maintains some of the nice properties of Haken manifolds. Most significantly, it still satisfies the geometrization conjecture. The idea of the proof is roughly as follows.  Let $M'$ be a Haken cover of a non-Haken manifold $M$.   The set of incompressible tori in $M'$ can be isotoped to be invariant under the covering transformations, so $M'$ cannot contain any incompressible tori since $M$ is non-Haken.  Therefore, $M'$ is either Seifert fibered or a hyperbolic manifold.   The case when $M'$ is Seifert fibered follows from \cite{ScottSFS}.  

If $M'$ is hyperbolic, then Mostow's rigidity theorem implies that the group of outer isomorphisms of $\pi_1(M')$ is isomorphic to a subgroup of the group of isometries of $M'$ \cite{thurstonnotes}.   We may assume that $M'$ is a regular cover of $M$ after passing to finite cover, which will still be hyperbolic. Then  any element $g$ of $\pi_1(M)$ not in $\pi_1(M')$ takes $\pi_1(M')$ to itself under conjugation, and hence is identified with an isometry of $M'$. Therefore there is a hyperbolic orbifold $N$ that has fundamental group $\pi_1(M)$.  $N$ is a actually a manifold  because  $N$ is a $K(\pi_1(M),1)$, and the fundamental group of such a complex has no torsion unless the complex is infinite dimensional. This is shown in \cite{Hempel} using the fact that the group homology of $\mathbb{Z}_n$, $H_i(\mathbb{Z}_n)$, is non-trivial for arbitrarily high $i$.  Since $M$ and $N$ are both covered by $\R^3$ they are both $K(\pi_1(M),1)$'s, and therefore are homotopy equivalent.   The main result of \cite{GMT} is that if $N$ is a closed hyperbolic 3-manifold, $M$ is an irreducible manifold, and $f: M \rightarrow N$ is a homotopy equivalence, then $f$ is homotopic to a homeomorphism.  Since we are in the situation of the hypothesis, $M$ is homeomorphic to a hyperbolic manifold when finitely covered by a hyperbolic manifold. 
 
 The virtual Haken conjecture asserts that if  $M$ is a closed irreducible three-manifold with infinite fundamental group,  then $M$ is finitely covered by a Haken manifold.  This implies that $M$ is either Seifert fibered or admits a hyperbolic structure. Much work has been done on this conjecture.  In particular, Dunfield and Thurston \cite{Nathanbill}  have shown that the 10,986 closed hyperbolic manifolds in the snappea census \cite{Snappea} are virtually Haken.  Cooper and Long proved that most Dehn fillings of a non-fibered  atorodial  three-manifold with torus boundary are virtually Haken in their paper, \cite{CooperLongvirt}.  Most recently, Masters, Menasco and Zhang \cite{MMZ} showed that infinitely many fillings of certain twist knots are virtually Haken. 
 
 In chapter 3, we show that infinitely many fillings of a large class of two-bridge knots are virtually Haken. 
The two-fold branched cover of the three-sphere, branched over a rational knot $\K{p/q}$, is the lens space $L_{p/q}$.  The preimage of the branching locus is an unknot $\LL$ in the lens space.  Furthermore, the $q$-fold cyclic cover of $L_{p/q}$ is the three-sphere  and $\LL$ lifts to a link $\D{p/q}$ in $S^3$.  The complement  $S_3 - \D{p/q} $ is thus the dihedral cover of $S^3 - \K{p/q}$. This cover is a great circle link complement.  We explain all this in more detail in chapter 3.  We also  show that when $p/q< 1/4$, the dihedral cover of the rational knot complement contains a closed incompressible surface $S$. We show that infinitely many fillings of the rational knot lift to fillings of the dihedral cover  where the incompressibility of this surface is preserved.  Thus, infinitely many fillings of a large class of  two-bridge knot complements are explicitly shown to be virtually Haken.  

Note that many rational knots  knots are fibered, and so do not fall under Cooper and Long's hypothesis.  
Specifically, it follows directly from the work of Gabai, \cite{detectingfiber} that a two-bridge knot $\K{p/q}$ is fibered if and only if $p/q = 1/(a_1+1/(a_2+1/(a_3 + .....a_n))) $ where each $a_i = \pm 2$.  Here we are using the fact that a rational knot $\K{p/q}$ with $p/q = 1/(a_1+1/(a_2+1/(a_3 + .....a_n))) $  can be written as the Murasugi sum of $n$ bands, where the ith band has $a_i$ half-twists, see \cite[Chapter 12]{Burde}. Using this, we can easily obtain infinitely many fibered two-bridge knots that are not torus knots and which can be written as $\K{p/q}$ with $p/q < 1/4$.  For example, the knot $\K{5/23}$ can be also be written as the knot $\K{18/23}$, using the equivalences of Lens spaces, and $18/23 = 1/(2+1/(-2+1/(2 +1/(-2 +1/(-2 +1/2)))))$.

%$18/23 = 1/(2+1/(-2+1/(2 +1/(-2 +1/(-2 +1/2)))))$.  Therefore, $\K{5/23}$ is a fibered knot such that infinitely many Dehn fillings of it are virtually Haken, since $5/23 <1/4$.  Similarly, $\K{7/33} \cong \K{26/33}$ and $26/33 = 1/(2 +1/( -2 +1/( 2+1/ -2 + 1/( -2 +1/( -2)))))$, while $7/33 <1/4$.  

\newpage
\pagestyle{myheadings}
\chapter{Classification}
\markright{  \rm \normalsize CHAPTER 2. \hspace{0.5cm}
  Classification }
\thispagestyle{myheadings} 
  
In this section, we introduce great circle links and exhibit some of their properties.  We characterize the torus decomposition of the complements of great circle links in $S^3$.  In addition, we classify great circle links up to five components. 

\section{Structure of the three-sphere} 

For the purposes of understanding the algebra of $S^3$, $S^3$ will be identified with the  unit quaternions, where the space of quaternions is  $\Q = \lbrace\qvec{q}= \quat{a}{b}{c}{d} : a,b,c,d \in \mathbb{R} \rbrace$ with $i^2 = j^2 = k^2 = -1, ij = k,jk = i,ki = j, ji = -k, kj = -i,ik = -j $, and the norm is given by $|\qvec{q} |^2 = a^2 +b^2 + c^2 + d^2$.  Any point ${\bf p}$ in the unit quaternions that has 0 real part is called a pure unit quaternion and we will call this 2-sphere \puq. Any such ${\bf p}$ defines two complex structures on the quaternions by left and right multiplication.  This is most commonly done with the pure unit quaternion $i$.   As such we have a notion of ${\bf p}$-complex lines in \Q $=\mathbb{R}^4$. $\lambda \qvec{x}$ and $ \qvec{x} \lambda$  will be called the left ${\bf p}$-line through $\qvec{x}$ and the right ${\bf p}$-line through $\qvec{x}$ respectively, where $ \lambda = l +m {\bf p}$ and $l, m \in \mathbb{R}$, and $\qvec{x}$ is any point in \Q.    When the context is clear we will omit the left or right.  The restriction of a left or right ${\bf p}$-line to the unit 3-sphere  $S^3$ will be called a left or right ${\bf p}$-fiber respectively.  These fibers are all geodesics.

\begin{lemma} Every geodesic $g$ in $S^3$ is a left fiber and a right fiber. 
\end{lemma}

%\begin{proof} 
$g$ is the intersection of a 2-plane in $\R^3$ with the unit 3-sphere in $\R^4$.  Any two points \qvec{x} and \qvec{y} on $g$ are linearly independent, so the plane can be defined as $\lbrace l\qvec{x} + m \qvec{y}: l ,m \in \mathbb{R} \rbrace$.  This plane intersects the two-sphere $\puq$ in at least two points. It also intersects the image of this sphere under left multiplication by \qvec{y} $\qvec{y}\puq$, in two points.  $\qvec{y}\puq$  is another two-sphere in $S^3$, since $\qvec{y} \in S^3$. Therefore $l\qvec{x} + m\qvec{y}=\qvec{yq}$ for some $\qvec{q} \in \puq$.  Furthermore, $l^2 +m^2 =1$, since  $ \qvec{x,y,q} \in S^3$. This implies that $\qvec{x}$ and \qvec{y} are in the right \qvec{q}-fiber $\qvec{y}\lambda$, where $\lambda= r+s\qvec{q}$ and $r^2 +s^2 =1$.  Note that $\qvec{y} =\qvec{y}(1 +0\qvec{q})$ and $\qvec{x} = \qvec{y}(-\frac{m}{l} +\frac{1}{l}  \qvec{q})$, so the right \qvec{q}-fiber through \qvec{y} is the exactly the geodesic of $S^3$ in the plane defined by \qvec{x} and \qvec{y}.  The identical method holds for writing $g$ as a left fiber.  $\Box$
%\end{proof}

 % In one of these fiber bundles the fibers wrap around this geodesic with a right-handed twist, and in the other fiber bundle the fibers wrap around with a left-handed twist. We will call these bundles right-handed and left-handed Hopf fibrations.  

%Multiplication of the right by ${\bf q} \in \puq$ takes a left ${\bf p}$-fiber to another left ${\bf p}$-fiber and visa versa. Left and right multiplication by elements of $\puq$ commute, $(\lambda_{\bf{p}}(\quat{a}{b}{c}{d}) ){\bf{q}} = \lambda_{\bf{p}}((\quat{a}{b}{c}{d} ){\bf{q}}))$, and the left $\bf{p}$-fiber of $\quat{a}{b}{c}{d}$ is taken to the left {\bf p}-fiber of $(\quat{a}{b}{c}{d}){\bf q}$ under multiplication on the right by \qvec{q}.
 
Consider the right \qvec{q}-fibers through \qvec{x} and \qvec{y}, $\qvec{x} (l + m \qvec{q})$ and $\qvec{y} (l + m \qvec{q})$, where $l^2 + m^2 =1$.   
The distance between two points \qvec{x} and \qvec{y} in $S^3$ is the change in angle along a geodesic that goes from one to the other.  Such a geodesic is a left ${\bf p}$-fiber for some  $\qvec{p} \in \Q$, so $(r + s \qvec{p}) \qvec{x} = \qvec{y}$.  Right multiplication by $(l+m \qvec{q})$ takes this left $\qvec{p}$-fiber to another \qvec{p}-fiber, preserving the distance between them, since  $(r + s \qvec{p}) \qvec{x} (l +m \qvec{q}) = \qvec{y} (l+m \qvec{q})$.  Therefore there is a well-defined distance between two right \qvec{q}-fibers, the shortest distance between a given point \qvec{x} on one fiber and any point on the other fiber. The same argument holds for left fibers. 

\begin{lemma} For every pair of points \qvec{p},\qvec{q} $\in \puq$, there is a geodesic in $S^3 $ that is a left $\qvec{p}$-fiber and a right \qvec{q}-fiber. 
\end{lemma}

%\begin{proof} 

First note  that the map $\psi: S^3 \rightarrow S0(4)$ by sending $\qvec{x}$ to the transformation of $S^3$ that sends \qvec{y} to $\qvec{xyx}^{-1}$ has image the natural copy of $SO(3) \subset SO(4)$ that fixes 1. 
This transformation leaves the 2-sphere of pure unit quaternions invariant, since it is a fixed distance from 1 in $\mathbb{R}^4$.   The kernal of $\psi$ is $\pm 1$, which is finite, so the image is 3-dimensional. Since $SO(3)$ is a connected 3-dimensional group, the image of $\psi$ is all of $S0(3)$.  The action of $SO(3)$ is transitive on the 2-sphere \puq, so there is a \qvec{x} $\in S^3$ such that $\qvec{p} = \qvec{xqx}^{-1}$.  We claim that the left \qvec{p}-fiber and the right \qvec{q}-fiber of \qvec{x} are the same.  The first is $\lbrace( l + m \qvec{p}) \qvec{x}| l,m \in \mathbb{R} ,l^2 + m^2 =1\rbrace$.  $( l + m \qvec{p}) \qvec{x} = l \qvec{x} + m \qvec{px} = l\qvec{x} +m \qvec{xq} = \qvec{x} (l +m \qvec{q})$, which is the right \qvec{q} fiber of \qvec{x}.  This proves the claim.  $\Box$

%\end{proof}

We take $\Q$ as a left \C-module.  Then $\Q = \C^2$ by $\phi: a + b\qvec{i} + c\qvec{j} + d\qvec{k} \mapsto (a + b\qvec{i}, c + d\qvec{i})$.  Stereographically project the unit three-sphere in $\C^2$ so that $(1,0)$ maps to the origin, and the circles $\lbrace (z,0), |z|=1 \rbrace$ and $\lbrace (0,w) \rbrace$ map to the $z$-axis $\cup \infty$ and the unit circle in the $x$-$y$ plane, respectively.  Then we can see that any two right $\qvec{p}$-fibers will wrap around each other with a left-handed twist and any two left $\qvec{p}$-fibers will wrap around each other with a right-handed twist. 

\begin{definition} The quotient of $S^3$ by identifying left (resp. right) \qvec{q}-fibers is $S^2$.  The induced fiber bundle structure will be called a {\it right-handed (resp. left-handed)  \qvec{q}-bundle}.  In this terminology, the standard Hopf fibration is a right-handed \qvec{i}-bundle. 
\end{definition}

Note that the distance between points on the quotient $S^2$ is the distance between \qvec{q}-fibers in $S^3$. 

\begin{lemma} 
\label{circlepoint}
The image of a geodesic of $S^3$ under the projection map associated to any right-handed or left-handed bundle is either a circle or a point. 
\end{lemma} 

%\begin{proof} 
Consider some geodesic $g$ in $S^3$.  Then $g$ is a left and a right fiber. Without loss of generality, consider a left \qvec{p}-bundle and look at the image of $g$ on the quotient 2-sphere. Since $g$ is a right \qvec{q}-fiber, $g = \lbrace \qvec{y}(l+m \qvec{q}) | l, m \in \mathbb{R}, l^2 +m^2 = 1 \rbrace$, for some \qvec{y} $\in S^3 $. $\qvec{y}$ is in some left \qvec{p}-fiber, $f_p$. The other fibers that $g$ intersects are  $f_p(l+m \qvec{q})$. (Recall that right multiplication by a unit quaternion takes left \qvec{p}-fibers to left \qvec{p}-fibers.) Let $h$ be the geodesic in $S^3$ that is both a right \qvec{q}-fiber and a left \qvec{p}-fiber.  Then any point in $g$ is a fixed distance $\alpha$ from $h$.  The  left \qvec{p}-fibers that $g$ intersects are also at distance $\alpha$ from $h$. Therefore, in the quotient $S^2$, the image of $g$ is a fixed distance from the image of $h$, a point on the sphere.  This can only be a point or a circle centered at $h$.  $\Box$ 

%\end{proof}

Note that a geodesic passes through each fiber of a generic geodesic fibration twice because if two geodesics intersect at a point \qvec{x}, then they also intersect at $\qvec{-x}$. 

\begin{definition}
We will call a link in $S^3$ such that each component is a geodesic of $S^3$ a {\it great circle link}.

\end{definition}

Note that any link $L$ in $S^3$ where each component is a round circle and the components pairwise link $\pm1$ can be realized as a great circle link, by the following argument due to Thurston.  Consider $S^3$ as the unit three-sphere in $\R^4$.  Then the components of the link $L$ are the intersection of planes in $\R^4$ with $S^3$.  Since each pair of components links once, the planes will intersect pairwise in a point in the interior of $B^3$, where $\partial B^3 = S^3$.  Now expand the unit $S^3$ to be a sphere with slightly larger radius.  The same planes in $\R^4$ will intersect this sphere in round circles that do not intersect, and pairwise link once.  Project these radially onto the unit three-sphere.  This will be a new link of round circles that pairwise link once that is isotopic to the original link.  These circles are the intersection of planes in $\R^4$ that intersect in the interior of the unit $B^3$ and their intersections are slightly closer to the origin than before.  Continue ``out to infinity'' until the planes do intersect exactly in the origin and we have isotoped our link to be a great circle link.  At no point in this process (including the limit) do the circles intersect.  Radial projection preserves the normalized distance between points.  Taking the intersection of a larger sphere with the same set of planes increases the normalized distance between points on different components, when these planes intersect in the interior of the unit $B^3$.

We will consider two great circle links to be equivalent if there is an isotopy between them through great circle links.  It is unknown if every two great circle links that are isotopic are isotopic through great circle links. However, we conjecture that there should be such a pair with a high number of link components.  The classification that follows shows that there is no such pair with five or fewer components. 

\section{Torus Decomposition} 
\begin{lemma} 
\label{unknotted} 
An incompressible torus in a great circle link complement is unknotted. 
\end{lemma} 

%\begin{proof}
We have a link of n great circles.   Pick some component.  Then the remaining links are contained in an open solid torus, and each travels once around in the longitudinal direction. Thus we may think of the link complement as the mapping torus $M_\phi$ of a homeomorphism $\phi: S \rightarrow S$, where  $S$ is a $p$ - punctured sphere and $\phi$ fixes each puncture.  Suppose there is an incompressible torus $T$ in $M$.    Then $T$ (in fact any incompressible surface)  can be isotoped so that the only singularities with respect to the fibers  are saddle singularities, by \cite{Hassmin}.  Then there is an induced foliation on the surface, the intersection of the fibration of $M_\phi$ by $p$-punctured spheres with the surface.  Saddle singularities will induce isolated singularities on the surface foliation with index -1.  By the Hopf index theorem,  the Euler characteristic of the surface is the sum of the indices of the isolated singularities of any codimension 1 foliation.  Therefore, saddles would imply additional genus, and $T$ must be vertical with respect to the fibration.   Pick some level, $F_0$.  Then $T \cap F_0$ is a union of circles.  An outermost circle must glue to an outermost circle under the monodromy $\phi$ since $\phi$ is an orientation preserving homeomorphism.   Any outermost circle must map to itself since $\phi$ takes each puncture to itself.  Therefore,  if there is one torus component, there will be one circle at each level.   At each level this circle bounds a disk on both sides, and the torus bounds a disk $\times S^1$ on both sides, and so is unknotted.   $\Box$ 
%\end{proof} 

We call a collection of non-boundary parallel incompressible tori $\it{minimal}$  if each component of the complement of the collection is either hyperbolic  or the union of fibers of the same Hopf fibration, and no subset of the collection would satisfy this condition. By lemma~\ref{unknotted}, these tori can be taken to be unknotted. Given a component of the complement of a minimal collection of tori, we can obtain a geodesic link in $S^3$ by adding a link component for each torus in the boundary.   We may assume that this new component is a great circle, since there were great circles on both sides of the torus before cutting.

Conversely, suppose that we have two geodesic links, $L_1$ and $L_2$ with $m$ and $n$ components. We can obtain a link with $m+n -2$ components in the following way.  Pick a link component $c_1$ of $L_1$ and a component $c_2$ of $L_2$.  Geodesic links in $S^3 $ can be regarded as the intersection of a plane in $\R^4$ with $S^3 = \lbrace (w,x,y,z )| w^2 +x^2 +y^2 +z^2 = 1 \rbrace$.  Without loss of generality, $c_1$ can be identified with the plane $(0,0, y,z)$. We can assume that there is a torus neighborhood $w^2 +x^2 \leq \delta$ around $c_1$ that does not intersect any of the other components of $L_1$. We move the other components by letting the planes associated with them flow under vector fields normal to the torus $w^2 + x^2 = \delta$. 

We claim that the circles associated with these planes will not intersect under this flow.   We will indicate a plane in $\R^4$ by a basis, $\lbrace (w_1,x_1,y_1,z_1), (w_2,x_2,y_2,z_2) \rbrace$. Two planes $P$ and $P'$ intersect in a 1-dimensional or greater vector space exactly when then their unit circles will intersect in a point or coincide.   At time $t$ , let plane $P$ and $P'$ be moved to the planes represented by $\lbrace (tw_1,tx_1,y_1,z_1), (tw_2,tx_2,y_2,z_2) \rbrace$ and $\lbrace (tw_1',tx_1',y_1',z_1'),(tw_2',tx_2',y_2',z_2') \rbrace$.   This  will not change the number of solutions to the homogeneous system:

 \[ \left [ 
\begin{array}{cccc}w_1 & w_2 & w_1' & w_2' \\x_1 & x_2 & x_1' & x_2' \\y_1 & y_2 & y_1' & y_2' \\z_1 & z_2 & z_1'& z_2'     
            \end{array} 
\right ] 
\left ( 
\begin{array}{l} a\\b\\c\\d \end{array}
\right ) =
\left (
\begin{array}{l} 0\\0\\0\\0\\ \end{array}
\right ) 
 \]

\begin{figure}[h]

\caption{ The torus sum of a left-handed Hopf link (-3) and a right-handed Hopf link (+3). 
\label{torussum}}
\centerline{\epsfig{file=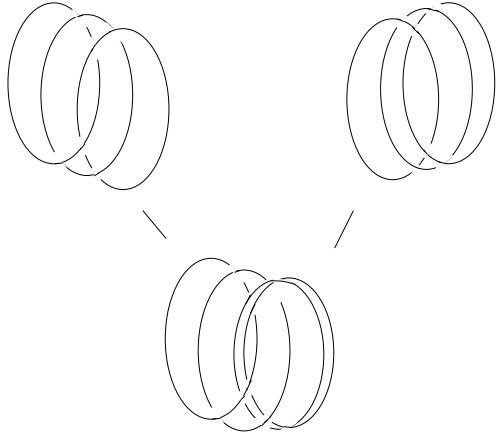, width=3in}}
\end{figure}

Eventually, the other planes intersect $S^3$ will satisfy $x_1^2 + x_2^2 \geq 1/\sqrt{2}$.  If we do the same thing for $L_2$, then we can glue the two boundary Clifford tori together to embed the link $\lbrace L_1 - c1  \rbrace \cup \lbrace L_2 -c2 \rbrace$ in $S^3$.  We will call such a link the {\it torus sum} of $L_1$ and $L_2$ along $c1$ and $c2$.  See figure \ref{torussum}. In this picture the bottom link is the torus sum of a left-handed Hopf link and a right-handed Hopf link along the 3rd component (from the left) of each link.  Swing the 3rd component of the right-handed link around to get the exact picture. 

\begin{lemma} Every geodesic link complement in $S^3$ can be written as a torus sum of pieces with hyperbolic or  Seifert fibered geometric structures.  The link components in each Seifert fibered piece are Hopf fibers of the same Hopf fibration, when considered as a link in $S^3$. 
\end{lemma} 

%\begin{proof} 

 By the Neilson-Thurston theory of surface automorphisms, $\phi$ is isotopic to an automorphism that is either reducible, periodic or pseudo-Anosov,  by ~\cite[Theorem 4]{thurstondiff}.  A pseudo-Anosov automorphism will yield a hyperbolic mapping torus, \cite{ThurstonHypII,Otal} and a periodic one will yield a Seifert fibered space. If $\phi$ is reducible, there are simple closed curves invariant under $\phi$ that divide the surface into regions where $\phi$ acts either as a periodic or pseudo-Anosov map.  These will result  in incompressible tori described above.  This proves the first claim of the lemma, which also follows from the geometrization of Haken manifolds, a proof of which can be found in \cite{Smithconj.}.

 Consider some piece that is Seifert fibered. Then the compact submanifold obtained by deleting an open regular neighborhood of each link component  is also Seifert-fibered.  This is a closed solid torus with some open neighborhoods of interior curves that intersect a meridian disk once removed.  The resulting local product structure can be extended  over the solid tori to give us a Seifert fibering of $S^3$. So now we need to consider Seifert fiberings of $S^3$, with base space $B$. If $B$ has a cover $\tilde B$, then any bundle structure over $B$ with fiber $F$ will extend to a bundle structure over $\tilde B$ with fiber $F$. This bundle will be a covering space of the bundle over $B$.  Therefore, since $S^3$ does not have any non-trivial covers, $B$ cannot have any non-trivial covers.  The base space $B$ is a two-dimensional orbifold and since it does not have any non-trivial covers it is either $S^2$, or a ``bad'' orbifold,  $S(p)$ or $S(p,q)$ where $p$ and $q$ are relatively prime \cite{orbifoldbook}. If $B = S^2$, then the Seifert fiber structure is a Hopf fibration and the link components are all fibers of this Hopf fibration.  If $B= S(p,q)$, then the regular fibers will link each other  with linking number $p \times q$ and the singular fibers with linking number $p$ or $q$.  Therefore, the only possible link components would be the singular fibers and there are only two, and any two geodesics in $S^3$ can be realized in their isotopy class as fibers of the same Hopf fibration.  The situation is the same  with $S(p)=S(p,1)$.  Therefore the Seifert fibered pieces all admit Hopf fibrations where the link components are fibers.  The associated  links are called Hopf links.   $\Box$

 %\end{proof}

From the above considerations, a great circle link is either a Hopf link, or a torus sum of Hopf links and/or hyperbolic great circle links, or a hyperbolic great circle link. 

\section{Classification}

Now we will classify great circle links with five or fewer components.  We know that pure positive and negative Hopf links occur, and that torus sums of lower component links occur.  We need to know when these various links are equivalent, and find all hyperbolic links.  Although we will use the orientation of $S^3$, we classify such links up to unoriented equivalence, and will change the orientation of a component when convenient.  

For the great circle links with two or three components, linear algebra is enough to classify them.  For four and five component links,  we will project a link along a Hopf fibration  and use the fact that all the link components project to round circles on the quotient two-sphere to classify the link. 

\subsection{Links with two and three components}

\begin{lemma} 

There is one two-component great circle link and there are two three-component great circle links. 
\end{lemma}

For two-component links: The two planes in $\mathbb{C}^2 = \R^4$ that define the link are linearly independent.  Therefore there is a linear transformation taking the two planes to the planes defined by $\lbrace (1,0,0,0),$ $ (0,1,0,0) \rbrace$ and $\lbrace (0,0,1,0),$  $(0,0,0,1) \rbrace$. 
  If this linear transformation is orientation preserving, it is isotopic to the identity through linear transformations. If it is not orientation preserving, we can flip the orientation of one of the components (by reordering the basis elements) so that it is orientation preserving. Topologically, this corresponds to the fact that if we have a link with two unknotted components that link $\pm1$, we can assume that they link $+1$ after possibly changing the orientation of one of the components. 

For three-component links:  Move the entire link by an isometry as above so that two of the links are in standard position.   The third component will be associated with some two-plane in $\R^4$ linearly independent from the other two.   This two-plane is the graph of a linear map from $\R^2 \mapsto \R^2$.   This is either orientation preserving or orientation reversing, and switching the order of the vectors that define the two-plane will not change this. 
Topologically, this represents the fact that the 3-component Hopf link  $6_3^3$ is chiral.  Therefore there are two distinct three-component great circle links, denoted by $+3$ and $-3$.  These are the two three-component links in figure \ref{torussum}.

\subsection{Links with four and five components}

For the case of four and five-component links, we will use the fact that the projection of a great circle link along a Hopf fibration is a collection of round circles on the quotient two-sphere, which follows directly from lemma \ref{circlepoint}.   There is a circle's worth of points in $S^3$ that map to a point on the two-sphere, and we will identify this circle with the unit tangent space at the point on the two-sphere.  In this way, we are identifying the geodesic with its quotient in $\mathbb{RP}^3$.  This is well-defined since every geodesic in $S^3$ is invariant under the antipodal map, and its quotient is a geodesic in $\mathbb{RP}^3$, the tangent space of $S^2$. In general the image of a component of a great circle link under the projection will be a round circle on $S^2$ with a vector field on it that is at a constant angle from the tangent vector field. We want to isotope through great circle links that do not intersect each other.  This exactly corresponds to making sure that no two vectors intersect.   

\begin{figure}[h]
\centerline{\epsfig{file=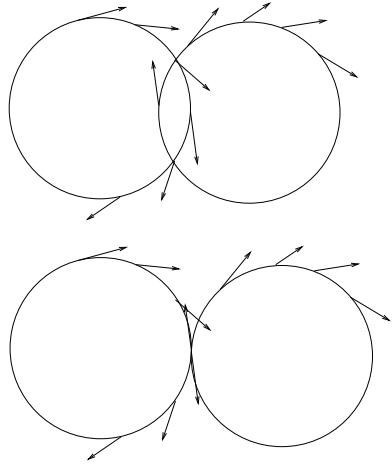,width=2in}}

\caption{The two great circles that map onto these circles can be isotoped to ``pull apart'' their projections.
\label{a}}
\end{figure}

\begin{figure}[h]
\centerline{\epsfig{file=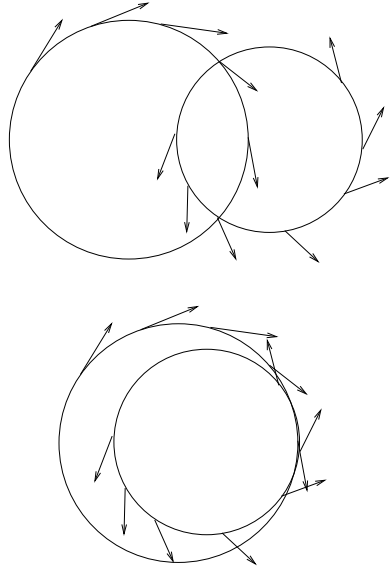,width=2in}}

\caption{The two great circles that map onto these circles can be isotoped so that the projection of one is inside the projection of the other.
\label{b}}

\end{figure}

Sometimes we do not need to know the vector field on the circle.  The image of a component of a great circle link on the two-sphere is a circle that divides the two-sphere into two components. If there are no other circles of the projection on one side,  we can isotope this component to a fiber without it intersecting the other components.  However, if the projections of two components intersect, there are two possibilities, distinguished by the vector fields on the projections. These two possibilities are pictured in figures \ref{a} and \ref{b}.

Note that the the first case we can push the two components apart from each other without their vector fields intersecting.  For the second one we can push one inside the other without their vector fields intersecting.  These two notions are the same on the two-sphere, but with other components in the projection, this will become important.
 
\begin{lemma} There are three four-component great circle links. 
\end{lemma} 

\begin{figure}[h]

\caption{ There are three four-component great circle links.
\label{class4}}
\input{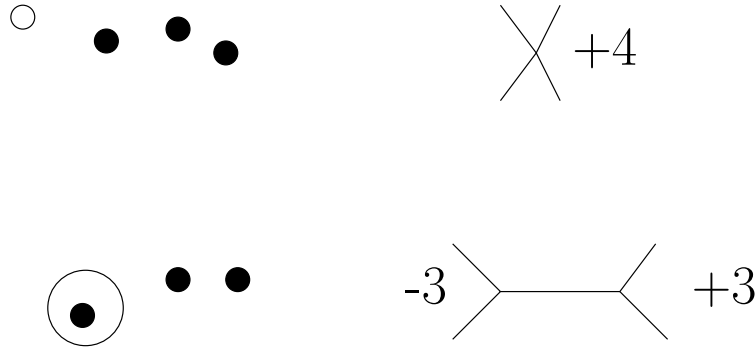}
\end{figure}

Let $L$ be a great circle link in $S^3$ with four components. Pick three of the components.  There is a linear transformation that takes these to a right-handed (positive) Hopf link or a left-handed (negative) Hopf link.  Transform the link by this linear transformation, and project by the Hopf fibration.  These three components of $L$ will project to points on the quotient two-sphere.  There is only one component of $L$ left, and it projects to a round circle.   This circle separates the two-sphere into two disks.  Either all three points corresponding to the other components are on one side of this circle, or the circle separates these points.   See figure \ref{class4} for the possible projections.

In the first case, no other components project onto one of the disks that this circle bounds.   There is a solid torus in $S^3$ that projects onto this disk under the Hopf fibration,  and we can isotope the component through great circles in this torus so that its projection is a point.  The link $L$ is a 4-component positive Hopf link if we projected by a right-handed (positive) Hopf fibration or a 4-component negative Hopf link, if we projected by a negative Hopf fibration.  

In the second case, concentric circles in a right-handed (left-handed) projection are the projections of left-handed (right-handed) Hopf links.  This can be seen directly by looking at the pre-image of such a configuration in  the solid torus above a disk that contains the concentric circles.  The vector fields on the circles determine the height mod $2 \pi$, which is well-defined in the solid torus because there is a section defined for the bundle everywhere except for a point.    So then this link is the torus sum of $+3$ and $-3$.  The torus is above a circle on the two-sphere that separates the concentric circles from the two points. As the order in which we do the torus sum does not matter,  there are three different four-component great circle links.

\begin{lemma} 
There are seven five-component great circle links. 
\end{lemma} 

\begin{figure}[h]

\caption{Projections of geodesic links with five components where at least four of the components are fibers of the same Hopf fibration.\label{atleast4}}
\input{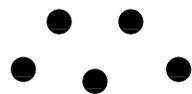}
\end{figure}

Again we will project along a Hopf fibration that contains three of the components.  Now there are two other components that project to circles on the two-sphere.  Assume first that one of these two components can also be isotoped to be a fiber so that there are four fibers.  Then there are three possibilities for the projection of the fifth component.  It either bounds a disk on one side that contains none of these four points or it separates one of the points from the other three, or it separates the four points into pairs.  These scenarios are pictured in figure \ref{atleast4}.  The first corresponds to the five-component positive  Hopf link, and the second to the torus sum of  a four-component positive Hopf link and a three-component negative Hopf link, assuming that the projection is by a positive fibration.  This is analogous to the situation for four components.  For the third configuration, note that there are two incompressible tori, one on the inside of the circle and another just on the outside.   If we start with three concentric circles, ($-3$) do torus sum on the component that projects to the outermost circle with $+3$, and also do torus sum on the component that projects to the innermost circle with $+3$, we will get a link that projects to the third configuration.  Therefore the third configuration corresponds to this torus sum when we project by a positive Hopf fibration.  

The same configurations for a projection by a negative (left-handed) Hopf fibration will correspond to analogous links with the signs of each piece in the torus sum reversed.

\begin{figure}[h]

\caption{Projections of geodesic links with five components where  three of the components are fibers of the same Hopf fibration.
\label{atmost3}}
\input{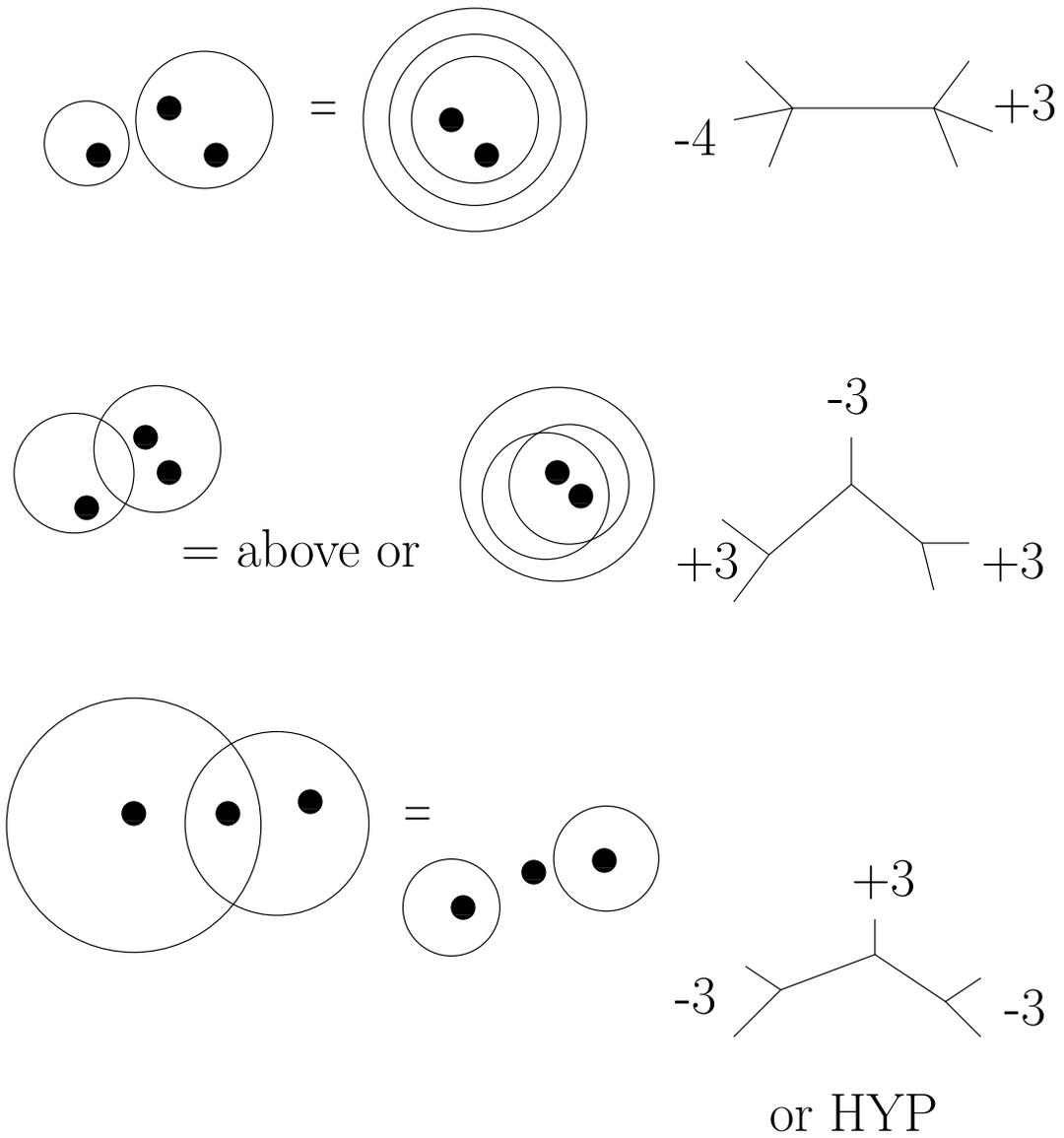}
\end{figure}

Now suppose that only 3 of the components project to points. Then there are two components that project to circles, and each of these circles separates the collection of three points on the two-sphere.  See figure \ref{atmost3} for the projections in this case.  If the two circles are disjoint, they still must bound disks containing points on both sides, and the two possibilities (the circles cobound an annulus or do not cobound an annulus) are at the top left and bottom right.  If the two circles intersect, they divide the sphere into four regions.  Then either one, two or three of these regions contain some or all of the three points.  If only one region contains points, then there are two adjacent regions that do not contain points, and the circles can be pushed one way or the other (see figures \ref{a} and \ref{b}) to become disjoint and then one will bound a disk on one side.  If two adjacent regions contain points, there will also be two empty adjacent regions.  The case when two non-adjacent regions contain points is the second row of figure \ref{atmost3} and the case when three regions contain points is the third row. 

 For the first row, after flipping one of the circles around the back of the two-sphere, we see the link is the torus sum of $-4$ and $+3$.  
 
For the second row, flip the component on the left around the back of the two-sphere.  Then we can either push this circle inside, as in figure \ref{b}, or not. If we can then we reduce to the case of the first row.  If we cannot, then we expand the isolated point to a small circle and flip it around.  Then we see that we have three concentric circles ($-3$) where we did torus sum along the components that project to the innermost circle and the second circle with $+3$.  

%Two fibers link positively in a positive Hopf fibration. In general, two great circles whose projections are intersecting circles link can be oriented to link positively if their projections be pulled apart and negatively if one of their projections can be pushed into the inside of the other (figures \ref{pushingoff1} and \ref{pushingoff2}, respectively). 

For the third row, we can flip the left circle around, push it into the interior of the second circle, and flip the second circle around to get $+3$ where we have performed torus sum with $-3$ along two components.  If we cannnot do the push, then there are no essential tori (recall incompressible tori are vertical) and our link is hyperbolic.  This link complement is the dihedral cover of the figure-8 knot complement, and will be discussed more later.

\newpage
\pagestyle{myheadings} 
\chapter{Some Large Great Circle Link Complements}
\markright{  \rm \normalsize CHAPTER 3. \hspace{0.5cm} 
  Some Large Great Circle Links}

\thispagestyle{myheadings}
In this section we introduce the great circle links $\D{p/q}$ for $q$ odd, $(p,q) = 1$.  We show that $S^3 - \D{p/q}$ is the dihedral cover of the rational knot complement $S^3 - \K{p/q}$, and that $S^3 - \D{p/q}$ is large when $p/q <1/4$.  We analyze the Dehn fillings of the manifolds $S^3 - \D{p/q}$ when $p/q \leq 1/4$. 
   
\section{The links  \D{p/q}}
\label{describelinks}
\begin{definition} The link $\D{p/q}$ in $S^3$ is the $q$ component great circle link  that is the orbit of the great circle $\lbrace (x + 0 i ,y + 0 i), x,y \in \mathbb{R}, x^2 + y^2 = 1 \rbrace$ under the action of the group generated by $\phi_{p/q}: (z,w) \mapsto (e^{{2 \pi i}/q} z, e^{{2 \pi i p}/q} w).  $  When it is clear from context, we will refer to $\D{p/q}$ by \D{} and $\phi_{p/q}$ as $\phi$. 

\end{definition}  

First note that these sets are the union of  great circles, since each component is the result of moving a great circle by an isometry.  We also claim that the images are disjoint. 

If not, then for some real  $r, r', s,s'$, $(e^{ \alpha i} r, e^{ \beta i} s) = (e^{ \alpha' i}  r', e^{\beta' i} s')$, where  $\alpha$ and $\alpha'$ are multiples of $\frac{2 \pi }{q}$ and $\beta $ and $\beta'$ are the respective multiples of $\frac{2 \pi  p}{q}$.  
This implies in particular that  $\tan \alpha = \tan \alpha'$, which only happens when  $\alpha = \alpha' \pm \frac{\pi}{2} \mod(\pi)$.  This never happens when $\alpha$ and $\alpha'$ are different multiples of $\frac{2 \pi}{q}$ with $q$ odd, the case we are considering. 
 
To understand these links, we will determine the image  of the 4 points $( \pm 1,0)$ and $(0, \pm1)$, under the action of $\phi$. The rest of each component will be the straight line (great circle) containing these points.

We will refer to a point  on one of the  ``core axes''   $\lbrace(z,w) \in S^3: w=0 \rbrace$ and $\lbrace(z,w) \in S^3: z=0 \rbrace$ by the angle the associated vector makes with ${(1,0)}$ or $(0,1)$ in the $z$ or $w$ complex plane in $\C^2$. 

\begin{figure}[htb]
\centerline{\epsfig{file=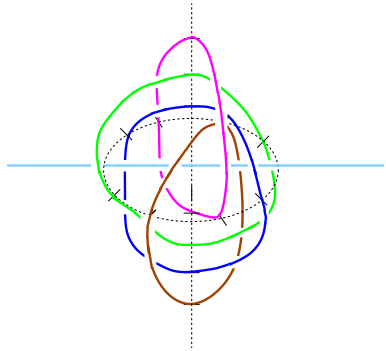, width=2in}}
\caption{The link $\D{2/5}$.  The $z$-axis of $S^3$ is the vertical line. 
\label{fig:5-2upright}}
\end{figure}

Consider the action of $\phi_{2/5}$ on the the real great circle, the circle that intersects the $z$-axis at $0$ and $\pi$ and the $w$-axis at $0$ and $\pi$. As we progress by $\phi_{2/5}$,  the components of the link $\D{p/q}$ will hit the  $z$-axis in the pairs of points $\lbrace(0,\pi) , ( \frac{2 \pi}{5},   \frac{7 \pi}{5}),( \frac{4 \pi}{5}, \frac{9 \pi}{5}),( \frac{6\pi}{5},  \frac{ \pi}{5}), (\frac{8\pi}{5}, \frac{3 \pi}{5}) \rbrace$.  The components intersect the $w$-axis in the pairs of points labeled $\lbrace(0,\pi), ( \frac{4 \pi}{5},   \frac{9 \pi}{5}),( \frac{8 \pi}{5}, \frac{3 \pi}{5}),$ $( \frac{2\pi}{5},  \frac{ 7\pi}{5}),$ $ (\frac{6\pi}{5}, \frac{ \pi}{5}) \rbrace$, in this order.  From this information, we can easily draw the link \D{2/5} as in figure \ref{fig:5-2upright}.  The $z$-axis is vertical and the first component is a horizontal line.

\begin{figure}[htb]
\centerline{\epsfig{file=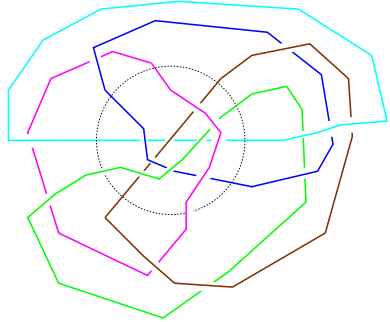, width=2in}}
\caption{A standard projection of the link $\D{2/5}$
\label{fig:5-2standard}}
\end{figure}

   Observe that  if $ n = \frac{q+1}{2}$, then $ n \frac{2 \pi}{q}  + \pi = \frac{\pi}{q} \mod 2 \pi$.  Therefore, as the intersection with the $z$-axis progresses by $ \frac{\pi}{q}$, the intersection with the $w$-axis progresses by $\frac{p \pi}{q}$. Number the components by the orbit of the real great circle $\lbrace (z,w) \in S^3: z \in \mathbb{R}, w \in \mathbb{R} \rbrace$ under the action of $\phi_{p/q}$.  The components numbered $0, \frac{q+1}{2}, 1, \frac{q+1}{2} +1, 2, \frac{q+1}{2} +2, ..... q-1 $ intersect the unit circle in the $z$ -plane in the points  labelled $0, \pi/q, 2\pi/q,3 \pi/q.... (q-1) \pi/q $ and their antipodal points.  The links in this order intersect the   the unit circle in the $w$-plane in the points labeled $0, p \pi/q, 2p \pi/q, 3p \pi/q.......(q-1) p \pi/q$, and their antipodal points.  

We will use this to draw a projection of these links in the plane. Place the unit circle in the $w$-plane on the $x-y$ plane and think of the unit circle in the $z$-plane as projecting to the $z$-axis in $\R^3$.  Then we draw the link with $q$ arcs connecting the pairs of points $(0, \pi), (\frac{ \pi p}{q},\frac{ \pi p}{q} + \pi), (\frac{ 2\pi p}{q},\frac{2 \pi p}{q} + \pi),  (\frac{ (q-1) \pi p}{q},\frac{(q-1) \pi p}{q} + \pi) $ inside the unit circle crossing consecutively crossing over all the previous arcs, and close them on the outside consecutively crossing under, so that they pairwise link once,  see figure \ref{fig:5-2standard}. This will be called a {\it standard projection} of the link \D{p/q}.

\section{Examples of surfaces in the link complements} 
We now order the components as they intersect the $w$-axis, starting with the point $(0, e^{0 i})$. The nth component will intersect the $w$-axis at a point that  makes an angle of  $n(\frac{\pi}{q})$ with $(0, e^{0 i})$. 
We define a genus $2p-1$ surface  $S_{n, n+2p-1}$
for each set of consecutive $2p$ components along the $w$-axis.   This surface consists of two $2p$-holed spheres, with $2p$ tubes connecting them.  We can also think of these surfaces as the boundary of a union of ``wedges'' and this will be useful later on. 

Define a {\it $z$-disk of radius $c$} to be a disk centered at a point $(e^{2 \pi i \theta},0)$  that consists of points $(z,w)$ in $S^3$ with $z$-coordinate $r e^{2 \pi i \theta}$, $r \in \R{}$ and $w$-coordinate satisfying $|w| \leq c$.  The boundary of this disk is a curve on the torus $\lbrace (z,w): |w| = c \rbrace = \lbrace (z,w) : |z| = \sqrt{1-c^2} \rbrace$.  A {\it $z$-wedge of radius $c$} is the union of $z$-disks for $z$ in some interval $(\alpha, \beta)$.  The center of a $z$-disk is a point in the $z$-axis, $(e^{2 \pi i \theta},0)$, and we will denote it by  $\theta$ when the context is clear.  $w$-disks and $w$-wedges are defined similarly.

Recall that a component of one of the links \D{p/q} looks like  $\lbrace (r e^{{n 2\pi i}/q}, s e^{{n 2\pi i p}/q})| r,s \in \mathbb{R}, r^2 +s^2 =1 \rbrace$.

\begin{lemma}
\label{zdiskfacts}
$z$-disk facts: 

\begin{enumerate} 

\item Two different $z$-disks (resp. $w$-disks) of radii $c$ and $c'$ intersect only when $c = c'= 1$.  Then their intersection is the $w$-axis (resp. $z$-axis).  

\item A $z$-disk of radius $c$ and a $w$-disk of radius $c'$ intersect only when $c^2 + c'^2 \geq 1$.  If $c^2 + c'^2 =1$, then they intersect in the point $(c e^{2 \pi i \beta}, c' e^{2 \pi i \gamma})$.

\item Every component in a dihedral link is in exactly two $z$-disks and two $w$-disks for any radii.
 For the component $\lbrace (r e^{{2\pi i}/q}, s e^{{2\pi i p}/q})| r,s \in \mathbb{R}, r^2 +s^2 =1 \rbrace$ these are the two $z$-disks centered at $\pm e^{{2\pi i}/q}$ and the two $w$-disks centered at   $\pm e^{{2\pi p i}/q}$.

\end{enumerate}

\end{lemma}
We can now easily describe the surface $S_{0, ... 2p-1}$.  It is the boundary of a handlebody with $2p$ open solid tori removed where the solid tori correspond to the components $\lbrace 0,...2p-1 \rbrace $.  Recall we are ordering the components as they intersect the $w$-axis, starting with the real great circle, and continuing counter clockwise. The surface $S_{0, ... 2p-1}$ is the boundary of the union of the $w$-wedges  of radius $c$ and $z$-wedges of radius $c'$, where $c^2 + c'^2 =1$.  The two $w$-wedges are the union of $z$-disks along  open intervals containing
$ [0 , \frac{(2p-1)\pi}{q}]$,  and  $[\pi, \frac{(2p-1)\pi}{q}+ \pi]$  where points on this axis are labelled by their angles.  We expand these intervals slightly so that the $w$-wedges contain arcs of the first $2p$ components in their interior.  Component $0$ intersects the $w$-axis in the points $(0, \pm e^{0 i })$ and component $2 p -1 $ intersects the $w$-axis in the points $(0, \pm e^{(2p -1)\pi i/q})$.  When $p/q < 1/2$, these intervals will not intersect. It follows from lemma \ref{covers}  and the classification of lens spaces that we can always express a dihedral link as $\D{p/q}$ with $p/q \leq 1/2$, so that this surface is always well-defined.  By the facts above, these two wedges do not intersect each other or any other component.  The two $w$-wedges are balls containing  exactly the first $2p-1$ components along the $w$-axis. 

Recall that components that intersect the $w$-axis in points $x$ and $x + \frac{p \pi}{q}$ and their antipodal points will intersect the $z$-axis in points $y$ and $y + \frac{\pi}{q}$ and their antipodal points. 
The $z$-wedges will be products of $z$-disks along intervals containing two consecutive components of the set $\lbrace 0,.... 2p-1 \rbrace$.  Every component in $\lbrace 0,... 2p-1 \rbrace $will intersect some such interval, since components $n$ and $n+p$ along the $w$-axis are next to each other along the $z$-axis.  We claim that there will be exactly $2 p$ $z$-wedges, each containing arcs of the two components that intersect $n$ and $n+p$ on the $w$-axis,  when $p/q <1/4$. The arcs labelled $n$ and $n+p$ are consecutive on the $z$-axis, so they are in the same $z$-wedge by construction.   These will be the only components in this $z$-wedge.  Indeed, the two components  immediately above and the two components  immediately below $n$ and $n+p$ along the $z$-axis are labelled $n+ 2p$, $n+3p$, $n-p$ and $n-2p$ and  these components are not in the set $\lbrace 0,... 2p-1 \rbrace$ if $n$ and $n+p$ are both in this set.  Thus there are exactly two components that intersect a given $z$-wedge defining the interior of $S_{0, ... 2p-1}$. Note that we have also shown that there are at least two components intersecting every complementary $z$-wedge.  In Figure \ref{fig:wedges} below, the two types of wedges will grow until they touch and the boundary of the resulting complex will be the surface $S_{0,1,2,3}$.  

\begin{figure}[htb]
\centerline{\epsfig{file=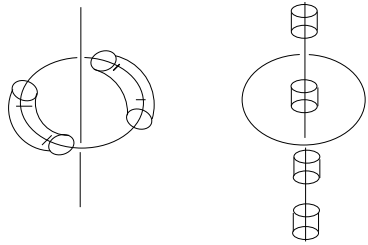, width=2in}}
\caption{$w$-wedges and $z$-wedges of the surface $S_{0,1,2,3}$
\label{fig:wedges}}
\end{figure}

\begin{figure}[htb]
\centerline{\epsfig{file=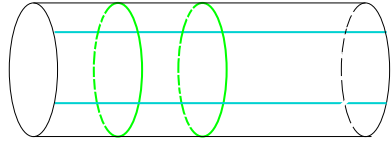, width=2in}}
\caption{Part of the torus where the two types of disks meet.  
\label{toruswithlines}}
\end{figure}

\begin{figure}[htb]
\centerline{\epsfig{file=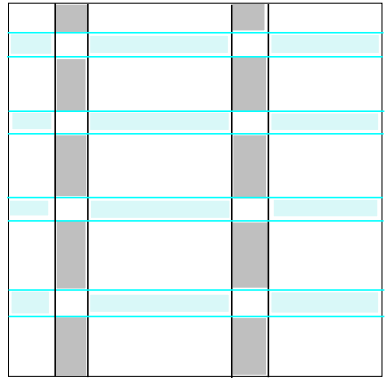, width=2in}}
\caption{Another view of the surface $S$. 
\label{surfaceontorus}}
\end{figure}

Although we have described the surface in terms of the handlebody in $S^3$  that it bounds, we can also describe it directly.  There are 4 $w$-disks, centered at just before $0$, just after $\frac{(2p-1) \pi}{q}$, just before $\pi$, and just after $\pi + \frac{(2p-1) \pi}{q}$.  There are $4 p$ $z$-disks that occur in pairs.  Each pair surrounds part of two components that intersect the $w$-axis at $x$ and $x + \frac{\pi p}{q}$.  These components intersect the $z$-axis consecutively. (Recall that as we progress by $\pi/q$ on the $z$-axis, we progress by $p \pi/q$ on the $w$-axis.)   The two different types of disks meet at the torus $\lbrace (z,w) \in S^3 : |z|=|w| \rbrace$.  The boundaries of the disks are  longitudes and meridians of this  torus, see figure \ref{toruswithlines}.  The surface on this torus is a checkerboard pattern connecting the two types of disks.   The meridianal curves bound disks of $S$ on the inside and the longitudinal curves  bound disks of $S$ on the other side of the torus.    The surface $S_{0, ..2p-1}$ will make a checkerboard pattern on this torus, as in figure \ref{surfaceontorus}.    The shadings mark the different sides of the surface.  The vertical arcs bound disks on one side, and the horizontal arcs bound disks on the other side.

\section{Incompressibility of the surfaces}

\begin{theorem} 
\label{incompressible}
$S_{0... 2p-1} $ is incompressible in $S^3 - \D{p/q}$ whenever $p/q < 1/4$. 
\end{theorem}

We will consider $M_{p/q}$ to be a the compact manifold  $S^3$ with  a regular neighborhood of each component of \D{p/q} removed.  $S_{0... 2p-1} $ is the closed surface of genus $2p-1$ as described above. 

First cut the manifold $M_{p/q}$ along  $S_{0... 2p-1} $ into two pieces, $M'$ and $M''$. If $S_{0... 2p-1} $ is incompressible in both $M' $ and $M''$,  then $S_{0... 2p-1} $  is incompressible in $M$.  

%For suppose there is a non-trivial curve on $S_{0... 2p-1} $ that bounds a disk $D$ in $M$.  Look at the intersection of this disk with $S=S_{0... 2p-1} $.  The intersection of $S$ and this disk is a collection of circles in $D$. The subdisk that any innermost circle bounds on $D$ is either a compressing disk for $S$, or $S$ can be isotoped to remove this circle of intersection.  Any innermost disk must be contained in $M'$ or $M''$.   So if there are no compressing disks in $M'$ or $M''$, then there are no compressing disks in $M_{p/q}$. 

Let $M'$ be the side where there are $2p$ components of the link. 

\begin{lemma}

When $p/q < 1/4$, $M''$  is homeomorphic to  $M'$ after trivially filling some of the torus boundary components of $M''$. 
\end{lemma}

 $M'$ and $M''$ are both the union of wedges.  The $w$-wedges of $M'$ are the two $w$-wedges that cover the points $\lbrace 0, \pi/q,2 \pi/q,...(2p-1)/q \rbrace$ their antipodal points. These are the points of intersection of the first $2p$ components along the $w$-axis and the $w$-axis. The two complementary wedges, the $w$-wedges of $M''$, contain at least $2p$ intersection points when $p/q <1/4$.  Similarly, there are $2p$ $z$-wedges of $M'$, each containing two arcs.  The  $2p$ complementary $z$-wedges  that make up the $z$-wedges of  $M''$ contain at least 2 arcs each when $p/q<1/4$.  In addition, the isometry $\phi_{p/q}: (z,w) \mapsto (e^{\frac{2 \pi i }{q}} z, e^{\frac{2 \pi i p}{q}} w) $ takes $M'$ into $M''$ when $p/q < 1/4$.    Note that every wedge of $M''$ contains exactly one wedge of $\phi(M')$.  Trivially fill in any torus boundary components of $M''$ that are not torus boundary components of  $\phi(M')$. Then this manifold is obtained from $\phi(M')$ by expanding the wedges of $\phi(M')$.  Thus $M'$ is homeomorphic to $M''$ with some of the torus boundary components of $M''$ trivially filled.  $\Box$

Suppose there was a compressing disk for $S$ in $M''$.  Since this disk does not bound a disk in $S= S_{0,...2p-1}$ and does not intersect any of the additional boundary components of $M''$, there is a compressing disk for $S$ in $M'$.   Therefore, it suffices to prove that there are no compressing disks in $M'$, the side with fewer components.  Note that $M'$ is the handlebody minus $2p$ open solid tori discussed in section \ref{describelinks}.  

Now we will cut $M'$ into two components, corresponding to the two $w$-wedges.  Each $z$-wedge meets a $w$-wedge in a two-holed disk.  The outer boundary component of this two-holed disk is a curve on the surface $S$.  For each $z$-wedge (which can be thought of as a tube connecting the $w$-wedges) cut along one of these two-holed disks.  After doing this to all the $z$-wedges, we will be left with 2 three-balls, $M'_+$ and $M'_-$,  minus neighborhoods of $2p$ arcs in each ball.   These correspond to the original two $w$-wedges. Each of $M'_+$ and $M'_-$ is decorated with the two-holed disks that we cut along.  Call the set of these two-holed disks $\mathcal{F}$.  Then $M'=  M'_+ \cup  M'_-$, glued along $\mathcal{F}$.

\begin{figure}[htb]
\centerline{\epsfig{file=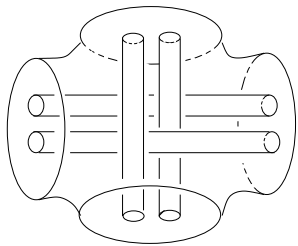, width=2in}}
\caption{Part of the surface $M'$ for $S_{0,1,2,3}$}
\label{fig:S_{0,1,2,3}}
\end{figure}

  Suppose that there is a compressing disk $D$ for $S$ in $M'$.   We will show:

\begin{enumerate}

\item A compressing disk $D$  can be isotoped off of $\mathcal{F}$, the set of 2-punctured disks that we cut along. 

\item There is no compressing disk in the complement of these  2-punctured disks.

\end{enumerate}

First we prove 1, by showing that the subset of the boundary consisting of these 2-holed disks is incompressible and boundary incompressible to $S$.

%\begin{figure}[htb]
%\centerline{\epsfig{file=twoholed.pdf, width=2in}}
%\label{fig: twoholed}
%\end{figure}

Let $F$ be any two-holed disk in $\mathcal{F}$, with one boundary curve meeting  $S$.  $F$ is incompressible in each piece of $M'$ cut along $\mathcal{F}$.  The fundamental groups of these pieces are $\ast_{2p} \Z$, the free product of $2p$ copies of $\Z$.  The map on fundamental groups induced by inclusion of $F$ into either of these pieces is the natural map $\Z \ast \Z \rightarrow \Z$ onto two components, which is an injective map.   Therefore, $F$ is incompressible in each piece of $M'$, since $F$ is two-sided, and each of $M'_+$ and $M'_-$ is orientable.  

Also, $F$ is boundary incompressible to $S$ in $M'$.  If not, then there is a non-trivial arc $c$ in $F$ isotopic through one of the pieces of $M'$ to an arc in $S$, where the endpoints of $c$ are fixed during the isotopy.  Note that the only non-trivial arc on $F$ that has boundary on $S$ is an arc going in between the two holes of $F$.  If this arc $c$ is isotopic to an arc $x$ on $S$, then $c$ together with an arc $s'$ on the boundary of $F$  is isotopic to $s' \cup x$, a loop on $S$.  This clearly cannot happen.  $c \cup s'$  links one tube once and all other tubes 0 times.  No such curve on $S$ does this, since any curve in $S$ is in the complement of the two-holed disks.  Thus $F$ is boundary incompressible in both pieces.

We will use the fact that each $F \in \mathcal{F}$ is incompressible and boundary incompressible in both parts of $M'$ to show that any compressing disk for $S$ can be isotoped off of $\mathcal{F}$.  Consider a compressing disk $D$ for $S$ in $M'$.  Then the intersection of $D$ and $\mathcal{F}$ is a collection of circles and arcs in  $D$.  Suppose that there is an innermost circle of this intersection in $D$.  Then this circle must bound a disk on some $F$, or else it would be a compressing circle.  Therefore, we can isotop $F$ to rid ourselves of this innermost disk, and continue until there are no circles of intersection. Now consider an outmost arc on $D$.  This arc is isotopic to an arc on $S$, namely part of the boundary of $D$.  Since $F \in \mathcal{F}$ is boundary incompressible,  it must be trivial in $F$, and we can isotop $D$ across this disk on $F$ to remove this arc of intersection.   Therefore, we can isotop a compressing disk $D$ to remove any outermost arc of  $\mathcal{F} \cap D$  in $D$ and hence all arcs. Thus we can isotop $D$ off all of $\mathcal{F}$.

 Recall that each two-holed disk corresponds to a $z$-wedge that intersects two  components of the link $\D{p/q}$.  Two such components will intersect the $w$-axis in points labeled $x$ and $x +\frac{p \pi}{q}$. Conversely, every two components of the link intersecting the $w$-axis in points labelled $x$ and $x + \frac{p  \pi}{q}$ will meet two of the two-holed disks, and these two disks will be called {\it corresponding two-holed disks} . 

To prove 2, suppose there is a compressing disk $D$ for  $S$ in $M'$ in the complement of $\mathcal{F}$. We cut along $\mathcal{F}$ and consider the piece that contains the disk.  The boundary curve $\partial D$  separates the collection of 2 punctured disks in this piece of $M'$. If all the two-holed disks are on one side, the curve is trivial.  If only one of the two-holed disks is on one side, then the disk $D$ separates corresponding two-holed disks, and intersects the associated components of the link \D{p/q}. Therefore the curve $\partial D$ must separate some corresponding pair of two-holed disks from some other corresponding pair of two-holed disks. These pairs intersect the $w$-axis in points labelled $x$ and $x + \frac{p \pi}{q}$ and $y$ and $y+  \frac{p \pi}{q}$, where $0 \leq x <y<x+ \frac{p \pi}{q}<y+ \frac{p \pi}{q} \leq  (2p-1)\frac{ \pi}{q} $.  This is the only possibility for the order since the wedge only covers the interval $[0,  (2p-1)\frac{p \pi}{q} ]$.  Ignoring the other components, the two pairs  in the $w$-wedge are positioned like those in \ref{fig:S_{0,1,2,3}}, with $x$ as the component closest to us, and $y+ \frac{p \pi}{q}$ as the component in the back. Then if there is a disk that separates the two pairs, we can double it along the boundary curve to get a sphere. This implies that there is a way to glue up the two corresponding pairs of two-holed disks so that the links are separated.  This is false, for once we glue one pair of corresponding disks, we will have a solid torus with two link components and two arcs in it. One of the remaining arcs is between the link components that are created, and cannot be isotoped relative to the boundary out of the solid torus. $\Box$

%Consider the manifold that is $S^3$ minus the links in $M'$ and the links in $\phi M'$.  These are the links that intersect the $w-$axis in the points labelled $\lbrace 0, \pi/q, 2\pi/q, .....(4p-1)\pi/q \rbrace$.  $S$ splits this manifold into two pieces homeomorphic to $M'$.  Any compressing disk for $S$ doubles in this manifold to be a disk that separates some of the link components from the others.  This is impossible because the links are all great circles, and must pairwise link once. Therefore there can be no compressing disk in $M'$ for $S$, and hence no compressing disk for $S$ in $S^3 - \D{p/q}$ when $p/q < 1/4$. $ \Box $ (theorem 3.3) 

\section{Dihedral covers of rational knots}
In this section, we will show that the link complements $S^3 - N(\D{p/q})$ are the  dihedral covers of the knot complements, $S^3 - N(\K{p/q})$.  

Note that the homotopy classes of arcs on a four-punctured sphere are indexed by their slopes.   Minimize the geometric intersection of an arc with the curves $m$ and $l$ below, where $m$ is the vertical slope, to put the curve in standard form.   We say an oriented arc has slope $p/q$ where $p$ is the oriented intersection with $l$ and $q$ is the oriented intersection with $m$ for some orientation of $l$ and $m$.   Without loss of generality, we can orient $m$ and $l$ so that a curve that twists in a right-hand fashion has positive slope.  This is similar to the case for curves on the boundary of a solid torus $T$ in $S^3$.  Let $\alpha$ be a curve that bounds a disk in $T$ and $\beta$ a curve that bounds a minimal genus surface in $S^3-T$ with boundary on $\partial T$.  Given some curve $c$ on  $\partial T$, the slope of $c$ is $p/q$, where $p$ is the oriented intersection with $\beta$  and $q$ is the oriented intersection with $\alpha$.  We can orient $\alpha$ and $\beta$ so that a curve of positive slope twists in a right-handed fashion.  Again, if we want to put the curve in standard form, we minimize the geometric intersection with $\alpha$ and $\beta$.   For the complement of the neighborhood of a link in $S^3$, we define the curves $\alpha$ and $\beta$  on each boundary component to be as if there were no other link components.  

\begin{figure}[htb]
\centerline{\epsfig{file=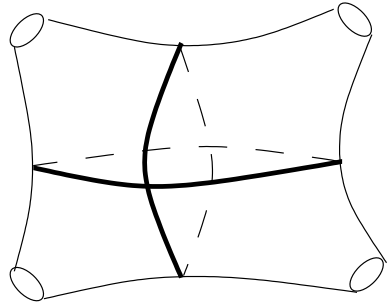, width=2in}}
\label{puncsphere}
\caption{The curves $m$ and $l$ on a four-punctured sphere.}
\end{figure}

\begin{definition}

Let $B$, resp. $B'$ be two balls with two standard interior vertical arcs marked. Let $h_{p/q}$ be the homeomorphism of their boundaries that takes a vertical arc (slope $\infty$) on $\partial B$ to an arc of slope $p/q$  on $\partial B'$.
The  rational knot $\K{p/q}$ is the union of the marked arcs in $S^3$ obtained by glueing together  $B$ and $B'$ by $h_{p/q}: \partial B \rightarrow \partial B'$. 
\end{definition}
\begin{figure}[htb]
\centerline{\epsfig{file=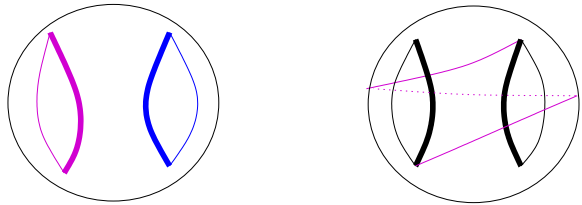, width=3in}}
\label{fig:knota}
\caption{Ex: the homeomorphism $h_{1/3}$ takes arc of slope $\infty$  to an arc of slope 1/3. }
\end{figure}

\begin{figure}[htb]
\centerline{\epsfig{file=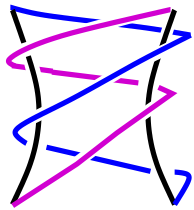, width=1 in}}
\label{fig:knotb}
\caption{The knot $\K{1/3}$, the trefoil.}
\end{figure}

\begin{lemma} 
Every two-bridge knot or link can be written as \K{p/q} for some $p/q, (p,q)=1$. 

\label{rationalform}
\end{lemma}

We will show that the two-fold cover of $S^2$ branched along a two-bridge knot is $L_{p/q}$. In the process, we will perform an isotopy on $S^3$ that takes a given two-bridge knot to rational form.  These two facts were first proven by Schubert, \cite{schubertknot}. 

Put the knot or link in two-bridge form in $S^3$.  Then the knot meets a 2-sphere $S$ in 4 points, $a,b,c,d$.  $S$ divides the 3-sphere into two balls each with two unknotted arcs of the knot in them.  We call these balls $B$ and $B'$, glued together by $h$ so that $h(\partial B) = \partial B'$, and $h(B) \cup B' = S^3$.    We may assume that the arcs in $B'$ (resp. $B$)  project to straight line segments (arcs of slope $1/0$) on $S$ ( resp. $h^{-1}(S)$.) 
Our goal is to show that $h$ is isotopic to $h_{p/q}$, for some $p,q, (p,q) = 1$.  $B$ and  $B'$ are two-fold  (branched) covered by solid tori $G$ and $G'$.  We'll identify each of the solid tori with $\lbrace (z,w ) \in S^3 : |z| <\frac{1}{\sqrt{2} } \rbrace$.  The group of deck transformations is $\Z_2$ and generated by  $\rho:(z,w) \rightarrow (\bar z, \bar w) $.  Note this fixes the great circle $\lbrace (z,w) \in S^3: z \in \mathbb{R}, w \in \mathbb{R} \rbrace$ and this intersects the solid torus in the two arcs  $A_1 = \lbrace (z,w) \in S^3: z \in \mathbb{R}, w \in \mathbb{R}, |z|<\frac{1}{2}, w >0 \rbrace$ and $A_2 =  \lbrace (z,w) \in S^3: z \in \mathbb{R}, w \in \mathbb{R}, |z|< \frac{1}{2}, w <0 \rbrace$.  The quotient of this solid torus by $\rho$ will be a ball with two standard arcs coming from $A_1$ and $A_2$.   Let  $\tilde h: \partial G \rightarrow \partial G'$ be the lift of $h$ that commutes with $\rho$,  $\rho \tilde h = h \rho$. 

\begin{figure}[htb]
\centerline{\epsfig{file=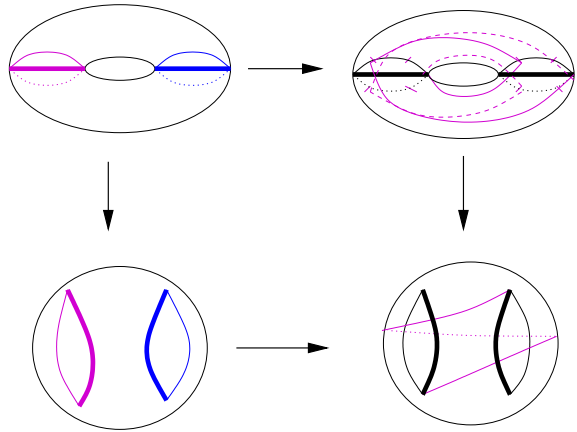, width=3in}}
\caption{The map $\tilde h$ commutes with the action of $\rho$}
\label{fig:knotlens}
\end{figure}

$\tilde h$ has the property that a meridianal curve on $\partial G$ wraps $p$ times around the z-direction (meridian) and q times around the w-direction (longitude) of $\partial G'$, for some $p,q \in \mathbb{Z}, (p,q) = 1$.   This $p/q$ determines the isotopy class of $\tilde h$.  However, we need this isotopy to commute with $\rho$.   Let $c_1$ be the meridianal curve $ \lbrace (\frac{1}{\sqrt{2}} e^{ i \theta}, \frac{1}{\sqrt{2}}) , 0 \leq \theta \leq 2 \pi \rbrace $ and $c_2 =  \lbrace (\frac{1}{\sqrt{2}} e^{ i \theta},- \frac{1}{\sqrt{2}}) , 0 \leq \theta \leq 2 \pi \rbrace $.  These will both map down to $1/0$ arcs on the sphere under the map $\rho: (e^{i \alpha}, e^{i \beta}) \mapsto (e^{-i \alpha}, e^{-i \beta})$.   First, minimize the intersection of $\tilde h(c_1(\theta))$ with the two meridianal curves $c'_1$ and $c'_2$ on $\partial G'$.   We say that two points $(x,y)$ on $\tilde h (c_1) \cap c'_1$ cancel if there is a disk bounded by the arc $xy$ on $\tilde h(c_1)$ and the arc $xy$ on $c'_1$, with no other points of $\tilde h(c_1)$ inside the disk.  If we have non-minimal intersection in $\tilde h (c_1) \cap c'_1$, then we will have such a disk, and since our map $\tilde h $ commutes with $\rho$, we have two such disks.  We can isotope $\tilde h$ in the union of  these two disks to reduce the intersection through an equivariant isotopy.  Do the same to reduce the intersections $\tilde h (c_1) \cap c'_2$, and $\tilde h (c_2) \cap c'_i, i = 1,2$.

The longitudinal curves $d'_i = \lbrace(\pm \frac{1}{\sqrt{2}}, \frac{1}{\sqrt{2}}e^{i \theta}) \rbrace$ on $\partial G'$ are also invariant under $\rho$. We can similarly reduce the intersection of $\tilde h(c_1)$, and then $\tilde h(c_2)$ with the $d_i$, through isotopies of $h$ that commute with $\rho$.  Note that the $c'_i$ and $d'_i$ divide $\partial G'$ into 4 quarters, paired by $\rho$.  We next completely straighten the image of the $c_1$ under $\tilde h$ to $\lbrace (e^{\frac{\theta p i}{q}}, e^{i \theta}), \theta \in [0, 2 \pi q) \rbrace$ and $c_2$ under $\tilde h$ to $\lbrace ( e^{\frac{\theta p i}{q}}, e^{i( \theta + \pi) }), \theta \in [0, 2 \pi q) \rbrace$, where $2 p$ is the minimal number of oriented intersections of $\tilde h(c_1)$ with $d'_1$, and $2q$ is the minimal number of oriented intersections of $\tilde h(c_1)$  with $c'_1$. $p$ is how many times the curve wraps around in the $z$-direction, and $q$ is  how many times  the curve wraps around in the $w$-direction.  We can make this isotopy that straightens the curve equivariant by copying the isotopy on one side of the $c'_i$  to the other. Now note that the complement of the two $(p,q)$ curves in $\partial G'$ is two regions that map to each other under $\rho$.  Also the pre-images of these two regions, the complement of the $c_i$ in $\partial G$, map to each other under $\rho$.  Therefore, we can straighten the rest of the map equivariantly by straightening it on one component of $\partial G - (\cup c_i) $ and extending the isotopy by $\rho$. 

Thus, the two-fold cover of $S^3$ branched along a two-bridge knot or link is a lens space, $L_{p/q}$. We will call such a two bridge knot  $\K{p/q}$, since the projection of a standard arc in one ball maps to a curve of slope $p/q$.   The knot complement, $S^3 - \K{p/q} $, is covered by $L_{p/q} - \LL$, where $\LL$      is the branching locus.  
$\Box$

\begin{lemma}
\label{covers}
Assume that $q$ is odd. 
$L_{p/q} - \LL$ is $q$-fold covered by the great circle link complement $S^3 -\D{p/q}$.
\end{lemma}

Recall that the link $\D{p/q}$ is the orbit of the great circle $\lbrace (z,w) \in S^3, z \in \mathbb{R}, w \in \mathbb{R} \rbrace $ under the action of $\Z_q$ generated by  $\phi_{p/q}: (z,w) \mapsto (e^{\frac{2 \pi i }{q}} z, e^{\frac{2 \pi i p}{q}} w )$. 
$S^3 - \D{p/q}$ modulo the action of $\phi_{p/q}$ is the Lens space $L_{p/q}$ minus a knot.  We claim that this knot is \LL.  A fundamental domain for this action is $F= (e^{i \alpha}, e^{i \beta}), 0 \leq \alpha, \beta < \frac{2 \pi}{q}$. The annulus $|z| = |w|$ in $F$ divides $F$ into a $z$-wedge and  a $w$-wedge, each of which glue up to a solid torus under the action of $\phi_{p/q}$. Each of these solid tori has two arcs removed, that are opposite from each other and contained in meridianal disks that come from $z$-disks and $w$-disks in the wedges.  $q$ vertical arcs on the boundary of the $z$-wedge will glue together to form one $p/q$ curve on the boundary of the resulting solid torus that connects the arcs on this side.   A meridianal curve on the boundary of the $w$-wedge will  glue to this $p/q$ curve.  
$\Box$

\begin{figure}[htb]
\centerline{\epsfig{file=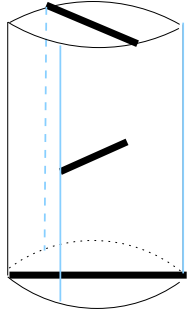, width=1in}}
\caption{This $z$-wedge will glue up to a solid torus with two arcs removed under $\phi_{1/3}$.  The blue curve will glue up to a (1,3) curve.}
\label{fig:knotlens}
\end{figure}

\begin{lemma}

Two two-bridge knots \K{p/q} and \K{p'/q'} are equivalent iff $q = q'$ and $p^{\pm1} = \pm p'(\mod q) $. 
\end{lemma}

We consider two knots in $S^3$ to be equivalent if there is a homeomorphism of $S^3$ that takes one knot to the other. 
Recall that the two-fold branched cover of $S^3$ branched along the knot $\K{p/q}$ is the lens space $L_{p/q}$.  
The conditions in the lemma are exactly the conditions for the Lens spaces $L_{p/q}$ and $L_{p'/q'}$ to be homeomorphic \cite{Brody}. If there is an isotopy of $S^3$ sending \K{p/q} to \K{p'/q'}, then the associated lens spaces are homeomorphic and $q = q'$ and $p^{\pm1} = p'(mod q)$.  On the other hand suppose $J: L_{p/q} \rightarrow L_{p'/q'}$ is a homeomorphism.  We will consider the lens space $L_{p/q}$ to be the union of two solid tori glued together by $\phi_{p/q}$ so that a meridian of one boundary torus glues to a $p/q$ curve on the other boundary torus. Then $J$ commutes with the action of $\rho$, since $\rho$ is in the center of $Map(T^2)$.  Therefore $J$ induces a homeomorphism $J' : L_{p/q} / \rho \rightarrow L_{p'/q'} / \rho$.  This is a homeomorphism of $S^3$ that takes one branching locus (the knot \K{p/q}) to the other branching locus (the knot \K{p'/q'}).  By definition, the two knots are equivalent.     $\Box$ 

\section{Virtually Haken Dehn fillings}

Next we discuss Dehn fillings of  rational knot complements $S^3 - N(\K{p/q})$.   We are interested in which of these fillings lift to fillings of the dihedral link complement.  Let $M_{p/q} = S^3 - N(\D{p/q})$, where $N(\D{p/q})$ denotes a regular neighborhood of the link.  We define the longitude of a component of $M_{p/q}$ to be a curve that bounds a disk in $S^3 $ minus that component. Note that $\phi_{p/q}$ takes a longitude of one component to a longitude of another. A longitude of $S^3 - N(\K{p/q})$  is a curve on the boundary torus that  lifts to a longitude of $L_{p/q} - N(\LL)$ which is a curve on the boundary torus that bounds a disk in $L_{p/q}$ which in turn lifts to the union of longitudes of $S^3 - N(\D{p/q}) $.   A meridian curve on the boundary of any of these manifolds is a curve that bounds a disk in the appropriate component of the complement. A curve that wraps $m$ times around the meridian of $M_{p/q}$ will lift  to a curve that wraps $1/2 m$ times around the meridian of $L_{p/q} -N( \LL) $ and this will lift to a curve that wraps $1/2 m$ times around each meridian curve of $S^3 -N( \D{p/q} )$.  Therefore, in order to lift a Dehn filling of the knot complement to a Dehn filling of the dihedral cover, $m$ must be even (where $m$ is the number of intersections of the filling curve with a longitude).  Specifically, a filling of the rational knot complement along a curve of slope $2m/l$ lifts to a filling of the dihedral cover along a multiple slope $\lbrace m/l, m/l, ...m/l \rbrace$.

We have shown that when a two bridge knot (or its mirror image) can be written as $\K{p/q}$ where $p/q \leq 1/4$, there is a closed incompressible surface $S$ in $M_{p/q}$, the dihedral cover  of $S^3 - N(\K{p/q})$.  We wish to determine when this surface remains incompressible under Dehn fillings.  Since $M_{p/q}$ has $q$  boundary components, a filling of $M_{p/q}$ is denoted by a multiple slope $\qvec{\gamma} = \lbrace \gamma_1, \gamma_2,... \gamma_q \rbrace$, where $\gamma_i$ is a slope on the ith torus boundary component $T_i$.  For each slope $\gamma_i =m/l$, we glue  a solid torus $F$ to $T_i$ so that a curve that bounds a disk on $\partial F$ is glued to a curve of slope $m/l$.  We use the convention that a filling corresponding to a curve of slope  $1/0$ is called trivial.   If $r$ and $s$ are slopes on a torus boundary component, $\Delta(r,s)$ denotes the minimum geometric intersection number between  two curves with slopes $r$ and $s$.  Let $\gamma$ be a multiple slope and $\beta =\lbrace \beta_1, ...\beta_n \rbrace$ a finite collection of imbedded non-trivial curves on $\partial M_{p/q}$.  Let 

$$\Delta(\gamma, \beta) = min_{i,j} \lbrace \Delta(\gamma_i, \beta_j), \beta_j \subset T_i \rbrace.$$

We can think of this as the distance between a complete filling and a partial filling in the case when there is at most one curve of $\beta$ on each boundary component.   Note that the collection $\beta$ could contain curves on any subset of the torus boundary components.

A slope on a boundary torus $T$ in a manifold $M$ that contains an incompressible  surface $S$ is said to be {\it coannular} to $S$ if there is an annulus $A$ in $M$ with one boundary component a non-trivial curve on $S$ and one boundary component on $T$.  The following theorem of Wu in \cite{Wu2} gives us sufficient conditions for an essential ($\pi_1$ - injective)  surface in a manifold with multiple torus boundary components to remain incompressible after Dehn filling. 

\begin{theorem}[Theorem 5.3 of \cite{Wu2}]
Let $X$ be a set of tori on the boundary of a compact, orientable hyperbolic 3-manifold $W$.  Let $F$ be a compact essential surface in $W$ with $\partial F \subset \partial W - X$, and let $\beta$ be the set of coannular slopes  to $F$ on $X$.  Then there is an integer $K$ and a finite set of slopes $\Lambda$ on $X$, such that $F$ is $\pi_1$-injective in $W(\gamma)$ for all multiple slopes $\gamma$ on $X$ such that $\Delta(\gamma, \beta) \geq K$ and $\gamma_i \notin \Lambda$.  

\end{theorem}

In light of this theorem, we have immediately:

\begin{theorem}
Infinitely many fillings of the rational knot complement $S^3 -N(\K{p/q})$ are virtually Haken when $p/q \leq 1/4$. 

\end{theorem}

The fact that infinitely many fillings are covered by manifolds that contain closed incompressible surfaces follows from Theorem \ref{incompressible} and the fact that the dihedral cover of a rational knot complement is  $M_{p/q} = S^3 - N(\D{p/q})$ (Lemma \ref{covers}).  Recall that any filling of a rational knot complement along a curve of slope $2m/l$ lifts to a filling of the dihedral cover along the multiple slope $\overline{m/l} = \lbrace m/l, m/l, ... m/l \rbrace$.  Then whenever $\Delta(\overline{m/q}, \beta) \geq K$  and $m/l \notin \Lambda$ where $\beta$ and $\Lambda$ are as above, then $S_{0,...2p-1}$ remains incompressible in $M_{p/q}(\overline{m/l})$.  According to \cite{hatcherthurston}, the only fillings of rational knots that do not produce irreducible manifolds are $\pm 2 q$ filling on $\K{p/q}$ with $p \equiv \pm 1 $(mod $q$) and $0$ filling on $\K{p/q}$ with $p \equiv 0 $(mod $q$).  A cover of an irreducible manifold is also irreducible.  The main result of \cite{hatcherthurston} shows that all but finitely many Dehn fillings of any rational knots are non-
Haken.  Therefore, inifinitely many fillings of a rational knot complement $\K{p/q}$ are non-Haken, and also covered by irreducible manifolds that contain a closed incompressible surface, when $p/q< 1/4$.  By definition, these fillings are virtually Haken.  $\Box$

We can make this more explicit by giving the set $\beta$ and the number $K$.   

\begin{lemma} 
Let $M_{p/q} = S^3 - N(\D{p/q})$ with $p/q<1/4$ and $S= S_{0,...2p-1}$. There are exactly four slopes on $\partial(M_{p/q})$ coannular to $S$, each of which has slope $\pm 1$. 
\end{lemma}

\begin{figure}[h]
\label{scheme}
\caption{How the link components intersect the $w$-axis}
\input{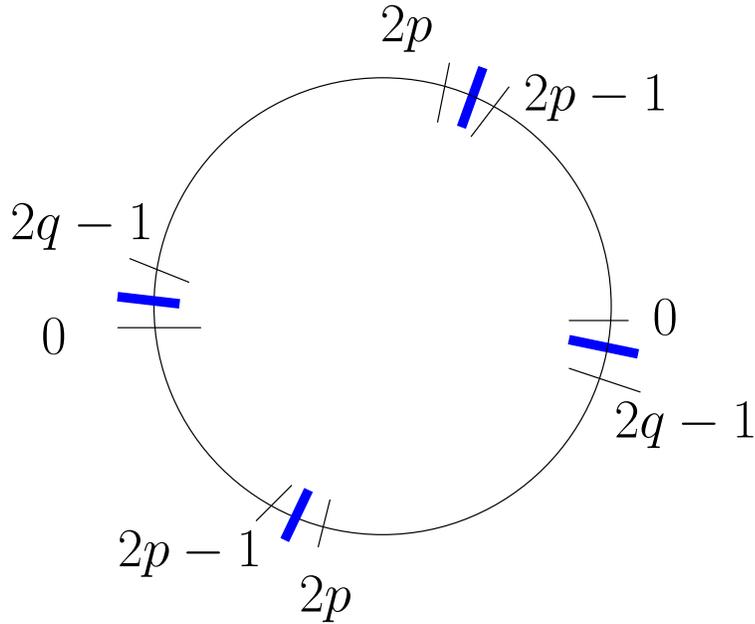}
\end{figure}

Order the boundary components as they occur along the $w$-axis as before.  Then components $0$ through $2p-1$ are all on the same side of $S$ in $M$.  Note that component $0$  is the same as component $2q$.  Components $0$, $2p-1$, $2p$ and $2q-1$ all are next to $w$-disks along the $w$-axis.  See the figure \ref{scheme} for how the surface $S_{0,...,2p-1}$  and the link components intersect the $w$-axis.  We claim that the boundary components corresponding to the link components $0$, $2p-1$, $2p$ and $2q-1$ are all coannular to $S$ and that  the coannular slope is $\pm 1$ in all these cases.  

Consider the components $0$ and $2p$.  All of the other components in their $z$-wedges intersect the $z$-axis above them, where we say  $(e^{\alpha_1}, 0) $ is above $(e^{\alpha_2}, 0) $ if $\alpha_1 > \alpha_2$ and the distance $|\alpha_1 - \alpha_2|$ is minimized among representatives of the angles.   We can see this from the standard projection construction. See figure \ref{fig:5-2standard} .  On the inside of the projection of the $w$-axis, components $x +p$, $x+2p$, $x+3p$, etc., cross over components $x$ and on the outside they cross under component $x$.  Therefore we can push components $0$ and $2p$ to the surface, by pushing them to the $w$-disks and $z$-disks simultaneously, without moving the other components.  The resulting annuli will each have a boundary component on a torus that goes once around the longitude and twists once in a left-handed fashion. 

Similarly,  components $2p-1$ and $2q-1$ are coannular to $S$.  All of the other components in their $z$-wedges are below them, since these components occur at the end of the $w$-wedge, after ordering the components counterclockwise.  They are isotopic to $S$ such that the resulting annuli will have boundary components on $\partial M_{p/q}$ that twist once in a right-handed fashion.

We claim that these are all of the coannular slopes.   Since we are ultimately interested in the incompressibility of the surface $S$ we need only consider annuli in a component of $M_{p/q} - S$.  Suppose that one of the other boundary components $c$ is coannular to $S$ in $M'$ or $M''$. Call the annulus $A$. First we reduce the intersection of $A$ and the collection $P$ of incompressible, boundary incompressible $n$- holed disks  that cut $M'$ or $M''$ into two balls with neighborhoods of arcs removed.    Trivial circles of $A \cap P$ on $A$ can be removed, starting with an innermost such circle, since $P$ is incompressible.  If there is a non-trivial  circle of $A \cap P$ on $P$, then this circle is also non-trivial in $A$, since $P$ is incompressible.  This implies this circle is boundary parallel on $A$, so we can isotop $A$ relative to the boundary to remove this circle of intersection. Now consider arcs of  $A \cup P$.  Call the 2 components of $P$ that touches $c$ $P_c$.  Then any arcs of $A \cap P$ are either connecting $c$ and $S$ on $P_c$ or trivial on $A$.  Since $P$ is boundary incompressible, we can isotope $A$ to remove these trivial arcs.  Then a component of the intersection of 
 $A$ with $M_{p/q} -P$ is a disk $D$.  $D$ gives us an isotopy of $c \cap  (M_{p/q} -P)$ into $S$ that does not move any other boundary components in $M_{p/q} - P$.

 Now if $c$ is not next to one of the $w$-disks along the $w$-axis, i.e., component $0$, $2p-1$, $2p$ or $2q-1$,  then there exists a component $x$ such that the components $\lbrace x, c,x+p,c+p \rbrace$ occur in that order (or $\lbrace c-p, x-p,c,x \rbrace$) along the $w$-axis, and these components are all in the same $w$-wedge.  Then we see that we are in the situation of figure \ref{fig:S_{0,1,2,3}}, where $c$ is the 2nd or 3rd component from the front.  Thus there can be no isotopy that moves $c$ into  $S$ fixing the other arcs. 

$\Box$

Since our surface is embedded, we can apply a theorem of Gordon and Luecke in \cite{CGLS}.  This says in part that if we have a closed incompressible surface in a manifold with a torus boundary component  with slope $r$ coannular to the surface, then Dehn filling along a slope $s$ with $\Delta(r,s) > 1$ leaves the surface incompressible.  Therefore, we can analyze a Dehn filling along a multiple slope as follows.   First fill in along the torus boundary components that do not have coannular slopes.  By Wu's theorem above, if we have avoided the finitely many slopes in $\Lambda$, our surface is still incompressible. Then we fill along the remaining torus boundary components which will still admit coannular slopes.  By the theorem 2.4.3 in \cite{CGLS}, the surface remains incompressible if we fill only along slopes that are at a distance greater than 1 from the coannular slopes.  Therefore, we may take $K$ to be 1.

%\newpage
%\pagestyle{myheadings} 
%\markright{  \rm \normalsize CHAPTER 4. \hspace{0.5cm}
% Some Small Great Circle Link Complements}
%\chapter{Some Small Great Circle Link Complements}
%\thispagestyle{myheadings} 
%\input{4_template.tex}

%%%%%%%%%%%%%%%%%    BIBLIOGRAPHY    %%%%%%%%%%%%%%%%%%%%%%%%%%%

%\newpage
%\pagestyle{myheadings} 
%\markright{  \rm \normalsize BIBLIOGRAPHY. \hspace{0.5cm}
  %YOUR TITLE}
%\begin{thebibliography}{99}
%\thispagestyle{myheadings}
%\addcontentsline{toc}{chapter}{\bf Bibliography}
%%%% here is a sample biliography entry %%%%%

%\bibitem{CE} H. ~Cartan and S. ~Eilenberg, {\it Homological Algebra}, Princeton University Press, Princeton (1956).

%\end{thebibliography}
\newpage
\pagestyle{plain}
\addcontentsline{toc}{chapter}{{\bf Bibliography }}
\bibliography{dissbib.bib}\bibliographystyle{alpha}

\begin{thebibliography}{CGLS87}

\bibitem[Bro60]{Brody}
E.~J. Brody.
\newblock The topological classification of lens spaces.
\newblock {\em Ann. of Math.}, 71:163--184, 1960.

\bibitem[Bur85]{Burde}
Gerhard Burde.
\newblock {\em Knots}.
\newblock Number~5 in De Gruyter studies in mathematics. W. De Gruyter, Berlin;
  New York, 1985.

\bibitem[CGLS87]{CGLS}
Marc Culler, C.~McA Gordon, J.~Luecke, and Peter~B. Shalen.
\newblock Dehn surgery on knots.
\newblock {\em Ann. of Math.}, 125:237--300, 1987.

\bibitem[CHK00]{orbifoldbook}
Daryl Cooper, Craig~D. Hodgson, and Steven~P. Kerckhoff.
\newblock {\em Three-dimensional Orbifolds and Cone-Manifolds}.
\newblock Number~5 in MSJ Menmoirs. Mathematical Society of Japan, Tokyo,
  Japan, 2000.

\bibitem[CHS96]{conwaypackings}
J.H. Conway, R.~H. Hardin, and N.~J.~A. Sloane.
\newblock Packing lines, planes, etc.: Packings in {G}rassmannian space.
\newblock {\em Experimental Mathematics}, 5:139--159, 1996.

\bibitem[CL99]{CooperLongvirt}
D.~Cooper and D.D Long.
\newblock Virtually {H}aken {D}ehn filling.
\newblock {\em J. Differential Geom.}, 52:173--187, 1999.

\bibitem[DT]{Nathanbill}
Nathan~M. Dunfield and William~P. Thurston.
\newblock The virtual {H}aken conjecture: experiments and examples.
\newblock arXiv:math.GT/0209214.

\bibitem[Gab86]{detectingfiber}
David Gabai.
\newblock Detecting fibred lnks in ${S}^3$.
\newblock {\em Commentarii Math. Helvetici}, 61:519--555, 1986.

\bibitem[GMT]{GMT}
David Gabai, G.~Robert Meyerhoff, and Nathaniel Thurston.
\newblock {Homotopy Hyperbolic 3-Manifolds are Hyperbolic}.
\newblock MSRI 1996-058.
\newblock arXiv:math.GT/9609207.

\bibitem[Has86]{Hassmin}
Joel Hass.
\newblock Minimal surfaces in foliated manifolds.
\newblock {\em Comment. Math. Helvetici}, 61:1--32, 1986.

\bibitem[Hem76]{Hempel}
John Hempel.
\newblock {\em 3-Manifolds}.
\newblock Number~86 in Annals of Mathematics Studies. Princeton University
  Press, Princeton, NJ, 1976.

\bibitem[HT85]{hatcherthurston}
A.~Hatcher and W.~Thurston.
\newblock Incompressible surfaces in 2-bridge knot complements.
\newblock {\em Inventiones Math.}, 79:225--246, 1985.

\bibitem[Lic62]{Lickorish}
W.~B.~R. Lickorish.
\newblock A representation of orientable combinatorial 3-manifolds.
\newblock {\em Ann. of Math.}, 76(3):531--540, 1962.

\bibitem[MMZ]{MMZ}
J.~Masters, W.~Menasco, and X.~Zhang.
\newblock Heegaard splittings and virtually {H}aken {D}ehn filling.
\newblock arXiv:math.GT/0210412.

\bibitem[Ota01]{Otal}
Jean-Pierre Otal.
\newblock The hyperbolization theorem for fibered 3-manifolds.
\newblock SMF/AMS Texts and Monographs, Vol. 7. American Mathematical Society/
  Soc\'ety Math\'ematique de France, 2001.

\bibitem[Sch56]{schubertknot}
H.~Schubert.
\newblock Knoten mit zwei br\"uken.
\newblock {\em Math. Zeitschr.}, 65:133--170, 1956.

\bibitem[Sco83a]{Scottsurvey}
Peter Scott.
\newblock The geometries of 3-manifolds.
\newblock {\em Bull. London. Math. Soc.}, 15:401--487, 1983.

\bibitem[Sco83b]{ScottSFS}
Peter Scott.
\newblock There are no fake {S}eifert fibre spaces with infinite $\pi_1$.
\newblock {\em Ann. of Math.}, 117:35--70, 1983.

\bibitem[Smi84]{Smithconj.}
In John~W. Morgan and Hyman Bass, editors, {\em The Smith Conjecture}, Pure and
  Applied Math., 112. Academic Press, Orlando, Fl., 1984.

\bibitem[Thua]{ThurstonHypII}
William~P. Thurston.
\newblock Hyperbolic structures on 3-manifolds, {II}: Surface groups and
  3-manfiolds that fiber over the circle.
\newblock arXiv:math.GT/9801045.

\bibitem[Thub]{thurstonnotes}
William~P. Thurston.
\newblock Princeton course notes, 1978-1979.
\newblock Available at http://www.msri.org.

\bibitem[Thu88]{thurstondiff}
William~P. Thurston.
\newblock On the geometry and dynamics of diffeomorphisms of surfaces.
\newblock {\em Bulletin of the American Mathematical Society}, 19(2):417--431,
  1988.

\bibitem[Thu97]{Thurstonbook}
William~P. Thurston.
\newblock {\em Three-Dimensional Geometry and Topology}.
\newblock Number~35 in Princeton Mathematical Series. Princeton University
  Press, Princeton, N.J., 1997.

\bibitem[VD89]{RP3}
O.~Ya. Viro and Yu.~V. Drobotukhina.
\newblock Configurations of skew-lines. (russian).
\newblock {\em Algebra i Analiz}, 1(4):222--246, 1989.
\newblock translation in Leningrad Math. J. 1 (1990), no. 4, 1027--1050.

\bibitem[Wal68]{waldhausenhaken}
Friedhelm Waldhausen.
\newblock On irreducible 3-manifolds which are sufficiently large.
\newblock {\em Ann. of Math.}, 87(1):56--88, 1968.

\bibitem[Wee]{Snappea}
J.~Weeks.
\newblock Snappea.
\newblock http://www.northnet.org/weeks/.

\bibitem[Wu]{Wu2}
Ying-Qing Wu.
\newblock Immersed surfaces and {D}ehn surgery.
\newblock {\em {\rm to appear in} Topology}.
\newblock arXiv:math.GT/9912049.

\end{thebibliography}

\end{document}